\documentclass[preprint]{elsarticle}

\usepackage[colorlinks]{hyperref}

\usepackage{amsmath,amssymb, amsfonts,graphics,xcolor,amsthm,dsfont}
\usepackage{epsfig}
\usepackage{latexsym}
\usepackage{euscript}
\usepackage{tocvsec2}
\usepackage{etoolbox, caption}
\usepackage{float}
\usepackage{marginnote}
\usepackage{tikz}
\usetikzlibrary{patterns}

\usepackage{enumitem}
\makeatletter
\newcommand{\mylabel}[2]{#2\def\@currentlabel{#2}\label{#1}}
\makeatother

\newtheorem*{theorem-no}{Theorem}

\usepackage{sectsty}
\allsectionsfont{\sffamily}
\usepackage{fancyhdr}
\newcommand{\comment}[1]{}  
\newcounter{pulse}[section]
\numberwithin{pulse}{section}  

\newcommand{\thf}{\sc} 

\theoremstyle{plain}
\newtheorem{theorem}[pulse]{\thf Theorem}

\newtheorem{proposition}[pulse]{\thf Proposition}
\newtheorem{lemma}[pulse]{\thf Lemma}
\newtheorem{claim}[pulse]{\thf Claim}
\newtheorem{corollary}[pulse]{\thf Corollary}
\theoremstyle{definition}
\newtheorem{definition}[pulse]{\thf Definition}
\newtheorem*{notation}{\thf Notation}

\theoremstyle{remark}
\newtheorem{remark}[pulse]{\thf Remark}

\def\UL{{\rm UL}}
\def\LR{{\rm LR}}

\def\B{\mathcal B}
\def\F{\mathcal F}
\def\G{\mathcal G}
\def\I{\mathcal I}
\def\R{\mathcal R}
\def\W{\mathcal W}
\def\cG{\mathbb G}

\def\bbN{\mathbb N}

\def\bbE{\mathbb E}
\def\bb1{\mathds 1}
\definecolor{Blue}{rgb}{0,0,1}
\definecolor{Green}{rgb}{0,0.9,0}
\definecolor{BlueGreen}{rgb}{0,0.5,1}
\definecolor{BrickRed}{rgb}{0.8,0.1,0.1}
\definecolor{RoyalPurple}{rgb}{0.6,0,0.6}
\definecolor{Brown}{rgb}{0.7,0.4,0.2}

\newcommand{\wep}{{\sl w.e.p. }}
\newcommand{\pink}[1]{{\color{black} #1}}

\newcommand{\blue}[1]{{\color{black} #1}}

\def\ch{}

\renewcommand\paragraph[1]{\par\vspace{1em}\noindent\textbf{#1}}

\begin{document}
\begin{frontmatter}

\title{Graph sequences sampled from Robinson graphons}

\author[mymainaddress]{Mahya Ghandehari\corref{mycorrespondingauthor}}
\cortext[mycorrespondingauthor]{Corresponding author}
\ead{mahya@udel.edu}
\author[mysecondaryaddress]{Jeannette Janssen}
\address[mymainaddress]{Department of Mathematical Sciences, University of Delaware, Newark, DE, USA, 19716}
\address[mysecondaryaddress]{Department of Mathematics \& Statistics, Dalhousie University, Halifax, Nova Scotia, Canada, B3H 3J5}

\begin{abstract}
The function $\Gamma$ on the space of graphons, introduced in \cite{CGHJK}, aims to measure the extent to which a graphon $w$ exhibits the {\sl Robinson} property: for all $x<y<z$, $w(x,z)\leq \min\{ w(x,y),w(y,z)\}$. Robinson graphons form a model for graphs with a natural line embedding so that most edges are local. \pink{The} function $\Gamma$ is compatible with the cut-norm $\|\cdot \|_\Box$, in the sense that graphons close in cut-norm have similar $\Gamma $-values. 
 \pink{In particular, any graphon close in cut-norm to the set of all Robinson graphons has small $\Gamma$-values.}
Here we show the converse, by proving that every graphon $w$ can be approximated by a Robinson graphon $R_w$ so that $\|w-R_w\|_\Box$ is bounded in terms of $\Gamma (w)$.
We then use classical techniques from functional analysis to show that a converging graph sequence $\{ G_n\}$ converges to a Robinson graphon if and only if $\Gamma (G_n)\rightarrow 0$.
Finally, using probabilistic techniques we show that the rate of convergence of $\Gamma$ for graph sequences sampled from a Robinson graphon can differ substantially depending on how strongly $w$ exhibits the Robinson property.
\end{abstract}

\begin{keyword}
graphon\sep Robinson \sep spatial network
\MSC[2020] 05C50 \sep  	05C62
\end{keyword}
\end{frontmatter}

\section{Introduction}
A graphon is a symmetric measurable function from $[0,1]^2$ to $[0,1]$. Graphons were introduced in \cite{LovaszSzegedy} as limit objects of convergent (dense) graph sequences. 
Every graphon $w$ gives rise to a rather general  random graph model 
 whose samples are graphs of any desired size.  
Such so-called $w$-random graphs are important as ``continuous'' generalizations of stochastic block models, which have been  prominent tools for modeling and analyzing large networks; see \cite{abbe} for a  survey on recent developments in the field of stochastic block models and their applications in community detection.  

In this paper, we focus our attention on graphons $w$ where the associated $w$-random graph model exhibits a spatial, line-embedded structure. In such models, vertices can be identified with points on the line segment $[0,1]$, and the link probability between points $x<y$ increases as $y$ moves closer to $x$.  
See for example \cite{Chang2019, Hoff2002} for such latent space models for social networks,  and \cite{smith} for an overview on continuous latent space models in general.
Graphons associated with such models are called \emph{Robinson graphons}. Namely, a  graphon $w:[0,1]^2\rightarrow [0,1]$ is {Robinson} if, for all $0\leq x\leq y\leq z\leq 1$, 
\begin{equation}\label{eq:Robinson}
 w(x,z)\leq \min\{w(x,y),w(y,z)\}.
\end{equation}
Robinson graphons were introduced in \cite{CGHJK} under the name of {diagonally increasing} graphons. 
The new terminology is compatible with the concept of Robinson matrices 
(also known as R-matrices) 
which appear in the study of the  well-known and challenging problem of \emph{seriation}. 
We refer the reader to \cite{Liiv2010} for a historical overview of seriation and its applications, 
and to \cite{Chepoi-seston11,Prea14,Fogel16,Laurent2017,FMR-statistical-seriation19} for recent advancements and new methodologies for seriation. 
Robinson matrices are significant from the graph theoretic point of view as well: the adjacency matrix of a graph, if labeled properly, is Robinson precisely when the graph is a unit interval graph. 

When dealing with real-life networks, an interesting question is whether a graph resembles an instance of a $w$-random process with Robinson $w$. 
One expects random instances of Robinson graphons  to be ``almost Robinson'' if the sample size is large. 
That is,  the sampled graph is expected to exhibit asymptotic \pink{aggregated} line-embedded behavior, in the sense that, \pink{for any two intervals $I_1\leq I_2$ in $[0,1]$, the density of edges joining vertices with labels in $I_1$ to those with labels in $I_2$ will be not significantly larger than the density of edges amongst vertices whose labels fall in the interval $[\sup(I_1),\inf(I_2)]$.}
Here, we pursue a systematic approach for formalizing these concepts.

This paper builds upon the previous work of the authors and collaborators in \cite{CGHJK}, where we introduced a function $\Gamma$ on the space of graphons such that 
\begin{enumerate}
\item\label{prop1} $\Gamma(w)=0$ precisely when $w$ is a Robinson graphon. 
\item\label{prop2} If $w$ is close to a Robinson graphon, then $\Gamma(w)$ is close to 0.
\end{enumerate}
Intuitively, we think of $\Gamma$ as a gauge of Robinson property. Here, the distance between graphons is measured by the cut-norm, denoted by $\|\cdot\|_\Box$, which has close relation with the theory of (dense) graph limits. Indeed, the cut-norm gives rise to the correct metric to define convergence of growing sequences of dense graphs in the sense of Lov\'{a}sz-Szegedy \cite{LovaszSzegedy}.
This is because the cut-norm is robust under sampling. Namely, for a growing sequence of graphs $\{G_n\}$ sampled from a graphon $w$, almost surely, the graphs in the sequence $\{G_n\}$ can be labeled  so that $\|w_{G_n}-w\|_\Box\rightarrow 0$.  Here $w_{G_n}$ is the graphon associated with the adjacency matrix of $G_n$; for details see Section~\ref{sec:def}.

Our original motivation for introducing $\Gamma$ in  \cite{CGHJK} was to recognize graph sequences which are sampled from Robinson graphons. To complete this task, we need to strengthen property~\ref{prop2} to an ``if and only if'' condition. This, however, proved to be a rather challenging problem, which we resolve in the current article.  
To our knowledge, $\Gamma$ is the first $\|\cdot\|_\Box$-norm continuous function which allows the recognition of graph sequences that are eventually Robinson. 
\pink{It may not be possible to compute the exact value of $\Gamma$, since it involves taking a supremum over all measurable sets.} However, 
the continuity results of \cite{CGHJK} guarantee that it can be approximated as accurately as one wants. 
This justifies our claim that $\Gamma$, or the graph/graphon parameter  it induces, is a suitable candidate for formalizing the notion of \emph{almost Robinson graphs}, i.e.~large graphs  sampled from Robinson graphons. 
%

\subsection{Contribution of this paper}
In this paper we show that (1) $\Gamma$ indicates closeness to a Robinson graphon, (2) $\Gamma$ 
identifies graph sequences converging to a Robinson graphon, and (3)  $\Gamma$-values of  graph sequences sampled from a Robinson graphon $w$ converge to zero at different rates depending on the strength of the Robinson property of $w$. 

\subsubsection*{I. Robinson approximation of graphons}
The main result of this article demonstrates the stability of $\Gamma$ near 0, in the sense that if $\Gamma(w)$ is close to 0, then $w$ must be close to a Robinson graphon. This resolves \cite[Conjecture~6.5]{CGHJK}.  Even though 
sections~\ref{sec:approx} and \ref{sec:proof-setting} deal with the particular function $\Gamma$,  the results hold for any $\Gamma$-type function, i.e. any function on the collection of graphons satisfying a natural condition  
(as in Lemma~\ref{lem:Gamma-lower-bound}), and attaining 0 when applied to Robinson graphons.
\begin{theorem}\label{thm:main-1}
Let $w:[0,1]^2\rightarrow [0,1]$ be a graphon. Then there exists a Robinson graphon $u$ so that
$\|w-u\|_\Box\leq 14\Gamma (w)^{1/7}$.
\end{theorem}
In Section~\ref{sec:approx},  we define a Robinson approximation $R_w$ for a graphon $w$ (see Definition~\ref{def:R(w)}), and in Section~\ref{sec:proof-setting}  we show that  $\|w-R_w\|_{\Box}\leq 14\Gamma (w)^{1/7}$. Even though the Robinson approximation  $R_w$ is easy to state, the proof is rather technical and involves delicate estimates. The significance of Theorem~\ref{thm:main-1} lies in the fact that it measures the error of approximation in terms of cut-norm, which is the suitable norm when studying  converging $w$-random graph sequences. The complicated nature of the proof is due to the facts that firstly  cut-norm is not as easy to compute as $\ell^p$ norms; and secondly, the function $\Gamma$ is defined via certain aggregated averages to ensure its continuity. 
Finally, suppose a given graph $G$ is  sampled from \emph{some} Robinson graphon. The Robinson approximation of $w_G$ provides an approximation for the underlying graphon. This is an instance of the \emph{graphon estimation problem}, where the goal is to invert the sampling process and to recover a graphon from a sampled graph.

\subsubsection*{II. Recognition of graph sequences sampled from Robinson graphons}
In Section~\ref{Sec:conjecture}, we present an important application of Theorem~\ref{thm:main-1}. Namely, we combine Theorem~\ref{thm:main-1} with some classical techniques from functional analysis to  prove the following:
\begin{theorem}\label{thm-app}
Let $\{G_n\}$ be a convergent sequence of (dense) graphs. If $\Gamma(G_n)\rightarrow 0$, then $\{G_n\}$ converges to a Robinson graphon. 
\end{theorem}
The limit object of a convergent graph  sequence  is an equivalence class of graphons defined by the cut-distance $\delta_\Box$, rather than the cut-norm itself (see Section~\ref{sec:def}). 
The cut-distance produces a  graph limit theory that applies to isomorphism classes of graphs.
In \cite{CGHJK}, it was shown that if  $\Gamma(G_n)\rightarrow 0$, then the limit object of $\{G_n\}$ can be represented by graphons  $u$ with arbitrary small $\Gamma(u)$. We were, however, unable to show that the limit object can be represented by an actual Robinson graphon; a gap which we close by Theorem~\ref{thm-app} of this article.
We show a similar result for the graph/graphon parameter induced from $\Gamma$, denoted by $\widetilde{\Gamma}$.
Theorem~\ref{thm-app}, together with property~\ref{prop2}, illustrates the significance of $\Gamma$ by proving that it identifies almost Robinson graphs.

\subsubsection*{III. Rate of decay for samples of a Robinson graphon}
Not all Robinson graphons exhibit the linearly embedded property to the same degree. Consider 
\begin{equation*}
\begin{array}{lcl}
w_1(x,y)=\left\{ \begin{array}{ll} 
p&\text{if }|y-x|\leq d,\\
0&\text{otherwise.}\\
\end{array}\right. 
&
\quad w_2(x,y)=
p-c|y-x|
\end{array}
\end{equation*}
where $p\in (0,1)$, $d\in (0,0.5)$,  and $c\in (0,p]$, 
which are both Robinson. Intuitively, however, $w_2$ has a stronger line-embedded representation. 
In Section~\ref{sec:decay}, we show that there is indeed a difference in order between `flat' graphons like $w_1$, and `steep' graphons like $w_2$.
\begin{theorem}\label{thm:convergence}
Let $w:[0,1]^2\rightarrow (0,1)$ be a Robinson graphon, and let $G\sim \cG (n,w)$. 
\begin{itemize}
\item[(i)] If $w$ has a \emph{flat} region, i.e.~there exist measurable sets $S,T\subseteq [0,1]$ with positive measure such that $w|_{\blue{S\times T}}=p$, then
with exponential probability for a graph $G$ sampled from $w$,  $\widetilde{\Gamma}(G)=\Omega ( n^{-1/2})$.
\item[(ii)] If $w$ is a \emph{steep} graphon, i.e.~its partial derivatives exist and are bounded away from 0, then with exponential probability
for a graph $G$ sampled from $w$,
$\widetilde{\Gamma}(G)= O ( n^{-2/3})$.
\end{itemize} 
\end{theorem}

\subsection{Similar graph parameters}
In \cite{GJseriation}, the authors introduced a function $\Gamma_1$ on the space of matrices, which attains 0 exactly when it is applied to a  Robinson matrix. While $\Gamma_1$ is easy to compute, it fails to be continuous in cut-norm (or equivalently the graph limit topology). So $\Gamma_1$ is not a suitable Robinson measurement for  growing networks,
 whereas $\Gamma$ provides us with a tool to measure Robinson resemblance of  large graphs.

Finally, we remark that this article tackles an instance of the question ``given a graphon with specific properties, how can we infer properties of the graphs which converge to it?'' These types of questions or their reverse versions have been studied for various classes of graphs/graphons; see for example \cite{bollobas2011} for graph sequences converging to monotone graphons and  \cite{DHJthreshold} for random threshold graphs.

\section{Definitions, notations and background}\label{sec:def}
We denote by \blue{$W_0$} the set of all measurable functions $w:[0,1]^2\rightarrow [0,1]$ which are symmetric, i.e. $w(x,y)=w(y,x)$ for every point $(x,y)$ in $[0,1]^2$.  Let $W$ denote the span of ${W}_0$,  i.e.~the set of all real-valued bounded symmetric and measurable functions on $[0,1]^2$. Functions in $\blue{{W}_0}$ are called graphons. Every $n\times n$ symmetric matrix $A=[a_{ij}]$ can be identified with a graphon, denoted by $w_A$, in the following manner: Partition $[0,1]$ into $n$ equal-sized intervals $I_i$. For every $i,j\in \{ 1,\dots ,n\}$, let $w_A$ attain $a_{ij}$ on $I_i\times I_j$. Every labeled graph $G$ on $n$ vertices can be identified with the graphon associated with the adjacency matrix of $G$. We denote this graphon by  $w_G$.


\subsection{Cut-norm, cut-distance, graph limits, and $w$-random graphs}
The topology described by convergent (dense) graph sequences can be formalized by endowing \blue{${W}$} with the cut-norm, introduced in \cite{FriezeKanan99}. For  $w \in \blue{W}$, the cut-norm is defined as:
\begin{equation}
\label{cutnorm} 
\| w \|_{\Box}= \sup_{S,T \subset [0,1]}\left|\int_{S \times T} w(x,y) dxdy \right|,
\end{equation}
 where the supremum is taken over all measurable subsets $S,T$ of $[0,1]$. To develop an unlabeled  graph limit theory, the cut-distance between $u,w\in\blue{W}$  is defined as follows.
 \begin{equation}
\label{cutdistance} 
\delta_{\Box}(u,w)= \inf_{\phi\in \Phi}\|u^\phi-w\|_\Box,
\end{equation}
where $\Phi$ is the space of all measure preserving bijections on $[0,1]$, and $w^\phi(x,y)=w(\phi(x),\phi(y))$. This definition ensures that $\delta_\Box(w,u)=0$ when the graphons $w$ and $u$ are associated with the same graph $G$ with two different vertex labelings. In general,  two graphons $u$ and $w$ are said to be \emph{$\delta_\Box$-equivalent} (or equivalent, for short), if $\delta_\Box(u,w)=0$.

It is known  that a graph sequence $\{ G_n\}$ converges in the sense of Lov\'{a}sz-Szegedy 
whenever the corresponding sequence of graphons $\{w_{G_n}\}$ is $\delta_\square$-Cauchy. The limit object for such a convergent sequence can be represented as a graphon in $\blue{W}_0$ (not necessarily integer-valued, or corresponding to a graph). That the graph sequence $\{G_n\}$ is convergent 
to a limit object $w \in \blue{W}_0$ is equivalent to $\delta_\Box(w_{G_n},w)\rightarrow 0$  as $n$ tends to infinity, which in turn is equivalent to the 
existence of suitable labelings for  vertices of $G_n$ for which we have
\begin{equation}
\label{cutdist} 
 \lVert w_{G_n} - w \rVert_{\square} = \underset{S,T \subset [0,1]} \sup \Bigl \lvert \int_{S \times T} ( w_{G_n} - w ) \Bigr \rvert \rightarrow 0.
\end{equation}
See \cite[Theorem 2.3]{BCLSV2011} for the above convergence results.

The concept of \emph{$w$-random graphs} was introduced in \cite{LovaszSzegedy}, as a tool for generating examples of convergent graph sequences. 
For a graphon $w$,  we define the random process $\cG(n,w)$ 
on the vertex set $\{ 1,2,\dots ,n\}$, where edges are formed according to $w$ in two steps. First, each vertex $i$ receives a label $x_i$ drawn uniformly at random from $[0,1]$. Next, for each pair of vertices $i<j$ independently, an edge $\{ i,j\}$ is added with probability $w(x_i,x_j)$. 
Such edge-independent $w$-random graphs arise naturally in the theory of graph limits. 
In fact, almost surely the sequence $\{\cG(n,w)\}_n$ forms a convergent graph sequence, for which the limit object is just the graphon $w$. See \cite{lovasz-book} for a comprehensive account of dense graph limit theory.

\subsection{Functions $\Gamma$ and $\widetilde{\Gamma}$}
We  now give the definition of  the function $\Gamma$, which is a non-negative valued function on $\blue{W}$. Note that $\Gamma$ is not a graphon parameter, as its value does not solely depend on the equivalence class of a given graphon, but rather on the actual representative itself. 
\begin{definition}[\cite{CGHJK}]
\label{def:Gamma}
For a function $w$ in $\blue{W}$, and a measurable subset $A\subseteq [0,1]$, we define
\begin{eqnarray*}
\Gamma (w,A)
&=&\int\int_{y<z}{\Big[} \int_{x\in A\cap [0,y]} \left( w(x,z)-w(x,y) \right)dx{\Big]}_+ dy dz\\
&+&\int\int_{y<z}{\Big[} \int_{x\in A\cap [z,1]}\left( w(x,y)-w(x,z)\right)dx{\big ]}_+ dy dz,\\
\end{eqnarray*}
where $[x]_+:=\max\{x,0\}$. 
Moreover, $\Gamma (w)$ is defined as 
$$\Gamma (w)=\sup\big\{\Gamma (w,A): A\subseteq [0,1] \mbox{ measurable}\big\}.$$ 
\end{definition}
The function $\Gamma$ attains 0 when applied to a Robinson graphon. 
It was shown in \cite[Proposition 4.2]{CGHJK} that  $\Gamma (w)=0$ if and only if $w$ is a.e.~equal to a Robinson graphon. Namely, $\Gamma (w)=0$ precisely when there exists a Robinson graphon $u$ so that $\|u-w\|_\Box =0$.  This fact is indeed a trivial case of 
Theorem~\ref{thm:main-1}.

One can think of $\Gamma$ as a function on labeled graphs in a natural manner, namely $\Gamma(G)=\Gamma(w_G)$ for a labeled graph $G$. 
When dealing with graphs, $\Gamma$ identifies unit interval graphs  labeled ``properly'', i.e.~labeled unit interval graphs  whose adjacency matrices are Robinson. 
We denote such labeling of a graph, if exists, a \emph{Robinson labeling.} 
To turn $\Gamma$ into a graphon/graph parameter, we consider the following natural definition:
 \begin{eqnarray}
 \widetilde{\Gamma}(w)&=&\inf\left\{\Gamma(u):\delta_\Box(u,w)=0\right\},\label{graphongammaT}\\
  \widetilde{\Gamma}(G)&=&\min\left\{\Gamma(w_H): \ H \mbox{ is a labeled graph isomorphic to } G\right\}.\label{graphgammaT}
 \end{eqnarray}
It was shown in \cite[Theorem 6.4]{CGHJK} that $\widetilde{\Gamma}$ is $\delta_\Box$-continuous. \blue{Note that $\delta_{\Box}(w_G,w_H)=0$
for any labeled graph $H$ isomorphic to $G$, so $\widetilde{\Gamma}(w_G)\leq \widetilde{\Gamma}(G)$. From Theorem 6.4 in \cite{CGHJK} we also have that, for any converging graph sequence $\{ G_n\}$, the numerical sequences $\{ \widetilde{\Gamma}(w_{G_n})\}$  and $\{\widetilde{\Gamma}(G_n)\}$ converge to the same value. We conjecture that, in general, $\widetilde{\Gamma}(G)\leq \widetilde{\Gamma}(w_G)+O(\frac1{\pink{|V(G)|}})$, 
but leave this question aside since it is tangential to the main theme of the paper. }

\begin{remark}\label{rem:diff-notation}
The concepts of $\Gamma(G)$ and $\widetilde{\Gamma}(G)$ were formulated  using a similar but discrete approach in \cite{CGHJK}, and were denoted by $\Gamma^*(G, \prec)$ and $\Gamma^*(G)$ respectively, where $\prec$  indicates the ordering on $V(G)$.
It turns out that the two approaches result in asymptotically equal functions; see Section~\ref{sec:decay} for more details. 
To declutter notations, we have chosen to work with the simpler definitions $\Gamma(G)$ and $\widetilde{\Gamma}(G)$ in the present article.
\end{remark}
\section{Robinson approximation of graphons}\label{sec:approx}
In this section, we define a Robinson approximation $R_w$ for a graphon $w$. Later, we will prove that the error of approximation,  namely $\|w-R_w\|_{\Box}$, is bounded in terms of $\Gamma(w)$. 
The following notations will be used in the rest of this article. 

\begin{notation}\label{notation-arerage}
Let $S$ and $T$ be measurable subsets of $[0,1]$. 
\begin{itemize}
\item[(i)] We write $S\leq T$  to signify that for all $x\in S$, $y\in T$, we have $x\leq y$.
\item[(ii)] The product $S\times T$ is called a \emph{cell} or \emph{rectangle} in $[0,1]^2$.  We defined the average value of $w$ on the rectangle $S\times T$ to be
$$\overline{w}(S,T)=\frac{1}{|S|\, |T|} \int_{S}\int_T w(x,y)\  \blue{dx\, dy}.$$
\end{itemize}
\end{notation}

We now proceed to define the Robinson approximation $R_w$ for a given graphon $w$. First, we need some preliminary definitions.  
Since both $w$ and $R_w$ are symmetric, we restrict our attention to the region above the diagonal, which we denote by $\Delta$. Namely, we denote
$$\Delta=\left\{(x,y)\in[0,1]^2:\ x\leq y\right\}.$$
The value of the Robinson approximation at a given point is determined by the behavior of the graphon on the {\sl upper left (UL)} and {\sl lower right (LR)} regions defined by that point. 
Precisely, for any  \blue{$(a,b)\in \Delta$}, we define
\begin{eqnarray*}
\UL(a,b)&=&[0,a]\times[b,1],\\
\LR(a,b)&=&[a,b]\times[a,b]\cap \Delta.\\
\end{eqnarray*}
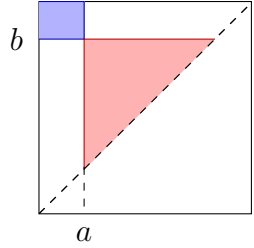
\begin{figure}[h]\label{fig:UL}
\begin{center}
\begin{tikzpicture}[scale=2]
\draw (0,0) -- (1.412,0) -- (1.412,-1.412) -- (0,-1.412) -- cycle;
\draw[dashed]  (0,-1.412) -- (1.42,0);
\filldraw[fill=blue!30, draw=blue!70!black] (0, 0)  -- (0.3,0) -- (0.3,-0.25) -- (0,-0.25);
\filldraw[fill=red!30, draw=red!70!black] (0.3, -1.12)  -- (0.3, -0.25) -- (1.17, -0.25);
\draw[dashed]  (0.3,-1.1) -- (0.3,-1.412);
\draw[dashed]  (0,-1.412) -- (1.42,0);
\node at (0.3,-1.55) {\large $a$};
\node at (-0.15,-0.25) {\large $b$};
\end{tikzpicture}
\end{center}
\caption{Regions $\UL (a,b)$ (blue)  and $\LR (a,b)$ (red)}
\end{figure}

The upper left and lower right regions provide us with an alternative, and rather more geometric, definition of Robinson graphons. 
Indeed, a graphon $w\in \blue{W}$ is Robinson if and only if, for all $(a,b)\in \Delta$ and all $(x,y)\in \UL (a,b)$, $w(x,y)\leq w(a,b)$. Alternatively, $w$ is Robinson if and only if, for all $(a,b)\in \Delta$ and all $(x,y)\in \LR (a,b)$, $w(x,y)\geq w(a,b)$. 

\begin{definition}[Robinson approximation for graphons]\label{def:R(w)} 
Given a graphon $w\in \blue{W}_0$ with $\Gamma (w) >0$, the \emph{Robinson approximation} $R_w$ is defined as follows. 
Let $\alpha=\Gamma(w)^{2/7}$.
Then for all $(x,y)\in \Delta$,
\begin{equation}
R_w(x,y)=R_w(y,x)=\sup\big\{ \overline{w}(S, T)\,:\, S\times T\subseteq \UL(x,y),\, |S|=|T|=\alpha\big\},
\end{equation}
taking the convention that $\sup\emptyset =0$. Moreover, we set $R_w=w$, if $w$ is Robinson itself.
\end{definition}
From the definition of UL and LR regions, it immediately follows that $R_w$ is indeed Robinson. Namely, let $(a,b)\in \Delta$ and $(x,y)\in \LR(a,b)$. Then $\UL(a,b)\subseteq \UL(x,y)$, and thus $R_w(a,b)\leq R_w(x,y)$.

We now restate our main theorem in a more detailed form. The proof will follow in Section~\ref{sec:proof-setting}.
\begin{theorem}[Equivalent to Theorem~\ref{thm:main-1}]\label{thm:main}
Let $w:[0,1]^2\rightarrow [0,1]$ be a graphon, 
and $R_w$ be as given in Definition~\ref{def:R(w)}. Then $\| R_w-w\|_{\Box}\leq 14\Gamma (w)^{1/7}$.
\end{theorem}
\section{Properties of $R_w$ and proof of Theorem~\ref{thm:main}}\label{sec:proof-setting}
The proof of Theorem~\ref{thm:main} relies on a  simple, yet important, lower estimate for $\Gamma(w)$, which we present in Subsection~\ref{subsec:Lemma} (see Lemma~\ref{lem:Gamma-lower-bound}). In fact, this lemma inspires the definition of the upper left (UL) and lower right (LR) regions. To obtain a bound on the error $\|w-R_w\|_\Box$, we analyze points $(x,y)\in \Delta$ based on the average  behavior of the graphon $w$ on the two regions 
$\UL(x,y)$ and $\LR(x,y)$. 
Indeed, we define black regions, white regions and grey regions in $\Delta$, each of which containing points with similar average behavior over the corresponding LR and/or UL regions. It turns out that  understanding features of these regions and their interactions  is the key to our final estimates; this is the content of Subsection~\ref{subse:region}.
The proof of Theorem~\ref{thm:main} is based on the idea that  the total area of all the grey regions is small (see Lemma~\ref{lem:grey-region-small}), while the local average difference between $w$ and $R_w$  inside either the black or white regions is controlled by $\Gamma(w)$. These facts lead to the conclusion that $\| R_w-w\|_\Box$ must be small.

\subsection{A lower bound on $\Gamma$ in terms of the ``upper left'' and ``lower right'' regions}\label{subsec:Lemma}
The lower bound for $\Gamma$, presented here, is in terms of average values of $w$ on rectangular regions of $[0,1]^2$ that are far from and close to the diagonal.  
\begin{lemma}\label{lem:Gamma-lower-bound}
Let $w\in \blue{W}$. Let  $S_u\leq S_l\leq T_l\leq T_u$ be subsets of $[0,1]$ with $|S_u|=|S_l|=|T_l|=|T_u| \geq \alpha$.
Then $$\Gamma(w)\geq \alpha^3 \Big(\overline{w}(S_u, T_u)-\overline{w}(S_l, T_l)\Big).$$
\end{lemma}

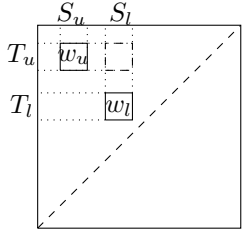
\begin{figure}
\centering
\begin{tikzpicture}[scale=0.3]
\draw (1,1) -- (10,1) -- (10,10) -- (1,10) -- (1,1);
\draw [dashed]  (1,1) -- (10,10);
\draw (4,5.8) -- (5.2,5.8) -- (5.2,7) -- (4,7) -- (4,5.8); 
\draw [dotted] (4,5.8) -- (1,5.8); 
\draw [dotted] (4,7) -- (1,7); 
\draw [dotted] (4,7) -- (4,10);  
\draw [dotted] (5.2,7) -- (5.2,10); 
\draw (2,8) -- (3.2,8) -- (3.2,9.2) -- (2,9.2) -- (2,8); 
\draw[dotted] (2,9) -- (2,10); 
\draw [dotted] (3.2,9) -- (3.2,10); 
\draw [dotted] (1,9.2) -- (5.2,9.2);    
\draw [dotted] (1,8) -- (5.2,8); 
\draw [dashed] (4,8) -- (5.2,8) -- (5.2,9.2) -- (4,9.2) -- (4,8); 
\node at (2.5,10.5) {$S_u$};
\node at (4.7,10.5) {$S_l$};
\node at (0.3,8.7) {$T_u$};
\node at (0.3,6.4) {$T_l$};
\node at (2.6,8.6) {$w_u$};
\node at (4.6,6.3)  {$w_l$};
\end{tikzpicture}
\caption{We let $w_u=\overline{w}(S_u\times T_u)$, $w_l=\overline{w}(S_l\times T_l)$. Lemma~\ref{lem:Gamma-lower-bound} bounds $\Gamma $ if $w_u>w_l$.}
\end{figure}

\begin{proof}
Setting $A=S_l\cup T_u$, we get
\begin{eqnarray*}
\Gamma(w)&\geq & \Gamma(w,A)\\
&\geq& \int_{x\in T_l}\int_{y\in T_u}\left[\int_{z\in S_l} w(y,z)-w(x,z)dz\right]_+dx\, dy\\
&+& \int_{x\in S_u}\int_{y\in S_l}\left[\int_{z\in T_u} w(x,z)-w(y,z)dz\right]_+dx\,  dy\\
&\geq&\left[\int_{x\in T_l}\int_{y\in T_u}\int_{z\in S_l} w(y,z)-w(x,z)dz\,  dx\,  dy\right]_+\\
&+& \left[\int_{x\in S_u}\int_{y\in S_l}\int_{z\in T_u} w(x,z)-w(y,z)dz\,  dx\,  dy\right]_+\\
&=&|T_l| |T_u| |S_l|\big[ \overline{w}(S_l,T_u)-\overline{w}(S_l, T_l)\big]_+ +
|S_l| |S_u| |T_u|\big[ \overline{w}(S_u, T_u)-\overline{w}(S_l, T_u)\big]_+\\
&\geq & 
\alpha^3\big[ \overline{w}(S_l,T_u)-\overline{w}(S_l, T_l)\big]_+ +
\alpha^3\big[ \overline{w}(S_u, T_u)-\overline{w}(S_l, T_u)\big]_+\\
&\geq&
\alpha^3 \big( \overline{w}(S_u,T_u)-\overline{w}(S_l, T_l) \big) .
\end{eqnarray*}
\end{proof}

\vspace{0.5cm}
\begin{remark}\label{remark:tight-example}
The following example shows that the above lower bound is sharp \blue{up to a factor of $\frac{5}{4}$}. For $i=1,\dots ,4$, let $I_i=(\frac{i-1}{4},\frac{i}4 ]$. Define $w\in \blue{W}_0$ to be
\[
w(x,y)=w(y,x)=
\left\{
\begin{array}{ll}
0 & \text{for } (x,y)\in I_2\times I_3,\\
\frac12 & \text{for } (x,y)\in I_1\times I_2\cup I_1\times I_4\cup I_3\times I_4, \\
\frac14 & \text{for } (x,y)\in I_1\times I_3\cup I_2\times I_4,\\
1 & \text{otherwise,}\\
\end{array}
\right.
\]
as shown in Figure~\ref{fig:exp-tight}.

\begin{figure}
\centering
\begin{tikzpicture}[scale=0.3]
\draw (1,1) -- (9,1) -- (9,9) -- (1,9) -- (1,1);
\draw [dashed]  (1,1) -- (9,9);
\draw [dotted] (1,3) -- (9,3);
\draw [dotted] (1,5) -- (9,5);
\draw [dotted] (1,7) -- (9,7);
\draw [dotted] (3,1) -- (3,9);
\draw [dotted] (5,1) -- (5,9);
\draw [dotted] (7,1) -- (7,9);
\node at (2,8) {\footnotesize{$1/2$}};
\node at (2,6) {\footnotesize{$1/4$}};
\node at (2,4)  {\footnotesize{$1/2$}};
\node at (2,2)  {\footnotesize{$1$}};
\node at (4,8)  {\footnotesize{$1/4$}};
\node at (4,6) {\footnotesize{$0$}};
\node at (4,4) {\footnotesize{$1$}};
\node at (4,2)  {\footnotesize{$1/2$}};
\node at (6,8)  {\footnotesize{$1/2$}};
\node at (6,6) {\footnotesize{$1$}};
\node at (6,4)  {\footnotesize{$0$}};
\node at (6,2)  {\footnotesize{$1/4$}};
\node at (8,8)  {\footnotesize{$1$}};
\node at (8,6)  {\footnotesize{$1/2$}};
\node at (8,4)  {\footnotesize{$1/4$}};
\node at (8,2)  {\footnotesize{$1/2$}};
\node at (0.2,2) {\footnotesize{$I_1$}};
\node at (0.2,4) {\footnotesize{$I_2$}};
\node at (0.2,6) {\footnotesize{$I_3$}};
\node at (0.2,8) {\footnotesize{$I_4$}};
\node at (2,0.2) {\footnotesize{$I_1$}};
\node at (4,0.2) {\footnotesize{$I_2$}};
\node at (6,0.2) {\footnotesize{$I_3$}};
\node at (8,0.2) {\footnotesize{$I_4$}};
\end{tikzpicture}
\caption{$w$ as in Remark~\ref{remark:tight-example}.}
\label{fig:exp-tight}
\end{figure}

Applying Lemma~\ref{lem:Gamma-lower-bound} to the sets $I_1\leq I_2\leq I_3\leq I_4$, we obtain that $\Gamma (w)\geq \left( \frac14\right)^3 (\frac12 )$. 
Now consider any set $A\subseteq [0,1]$, and let $A_i= A\cap I_i$ and $\alpha_i=|A_i|$ for $i=1,\dots ,4$.
If $y,z\in I_i$ for some $i=1,\dots ,4$, then $w(x,z)=w(x,y)$ for all $x$.
If $y\in I_1\cup I_2$, then $w(x,z)\leq w(x,y)$ for all $x<y<z$. If $z\in I_3\cup I_4$, then $w(x,y)\leq w(x,z)$ for all $x> z>y$. Therefore, only pairs $(y,z)\in I_3\times I_4\cup I_1\times I_2$ can make a positive contribution to $\Gamma (w,A)$, and
\begin{eqnarray*}
\Gamma (w,A) &=& \int_{y\in I_3}\int_{z\in I_4}\Big[ \int_{x\in A_1} (\frac12-\frac14)dx+ \int_{x\in A_2} (\frac14-0)dx
+ \int_{x\in A_3\cap [0,y]} (\frac12-1)dx\Big]_+\\
&& + \int_{y\in I_1}\int_{z\in I_2}\Big[ \int_{x\in A_2\cap [z,1]} (\frac12-1)dx+ \int_{x\in A_3} (\frac14-0)dx
+ \int_{x\in A_4} (\frac12-\frac14 )dx\Big]_+ \\
 &=&\frac{1}{16} \left(\int_{y\in I_3}\Big[\alpha_1+\alpha_2- 2\big|A_3\cap [\frac{1}{2},y]\big|\Big]_+dy
+\int_{z\in I_2}\Big[\alpha_3+\alpha_4-2\big|A_2\cap [z,\frac{1}{2}]\big|\Big]_+dz\right)\\
&\leq&\frac{1}{16}\left(\int_{y=\frac12}^{\frac34}\Big[\alpha_1+\alpha_2-2[y+\alpha_3-\frac34]_+\Big]_+dy+ \int_{z=\frac14}^{\frac12}\Big[\alpha_3+\alpha_4-2[\frac14+\alpha_2-z]_+\Big]_+dz\right)\\
&=&\frac{1}{16}\left(\blue{\int_{y=\frac12}^{\frac34-\alpha_3}(\alpha_1+\alpha_2)\ dy+}\int_{y=\frac34-\alpha_3}^{\frac34}\Big[\alpha_1+\alpha_2-2(y+\alpha_3-\frac34)\Big]_+dy\right.\\ 
&& + \left. \int_{z=\frac14}^{\frac14+\alpha_2}\Big[\alpha_3+\alpha_4-2(\frac14+\alpha_2-z)\Big]_+dz + 
\blue{\int_{\frac14+\alpha_2}^{\frac{1}{2}}(\alpha_3+\alpha_4)dz}  \right)\\
&=&\frac{1}{16} \left(\int_{y=0}^{\alpha_3}\Big[\alpha_1+\alpha_2-2y\Big]_+dy+ \int_{z=0}^{\alpha_2}\Big[\alpha_3+\alpha_4-2z\Big]_+dz\right.\\
&& \left. \blue{+(\alpha_1+\alpha_2)(\frac14-\alpha_3)+(\alpha_3+\alpha_4)(\frac14-\alpha_2)}\right),
\end{eqnarray*}
where we used the fact that, for any $y\in I_3$, the measure $\big|A_3\cap [\frac12,y]\big|$ is minimized when $A_3=[\frac34-\alpha_3,\frac34]$, in which case we have
$\big|A_3\cap [\frac12,y]\big|= [y-\frac34+\alpha_3]_+$ for every $y\in I_3$.
Similarly, for $z\in I_2$,
$\big|A_2\cap [z,\frac{1}{2}]\big|$ is minimized when $A_2= (\frac14,\frac14+\alpha_2]$, and the minimum value equals to $[\frac14+\alpha_2-z]_+$.
\blue{The right hand side of the above inequality increases when  $\alpha_1=\alpha_4=\frac14$. So we have 
\begin{eqnarray}\label{eq:ex-bound}
\Gamma (w,A)&\leq&\frac{1}{16} \left(\int_{y=0}^{\alpha_3}\Big[\frac14+\alpha_2-2y\Big]_+dy+ \int_{z=0}^{\alpha_2}\Big[\frac14+\alpha_3-2z\Big]_+dz
+\frac18-2\alpha_2\alpha_3\right).
\end{eqnarray}
By symmetry, we may assume \emph{wlog} that $\alpha_2\leq \alpha_3$. Note also that $\alpha_3\leq \frac18 +\frac{\alpha_3}2$, since $\alpha_3\leq \frac14$. In the case where $\alpha_3\leq \frac18+\frac{\alpha_2}{2}$, the right hand side of \eqref{eq:ex-bound} reduces to 
\[
\frac{1}{16}\left((\frac14+\alpha_2)\alpha_3 -\alpha_3^2 + (\frac14+\alpha_3)\alpha_2-\alpha_2^2+\frac18-2\alpha_2\alpha_3\right).
\]
This is maximized when $\alpha_2=\alpha_3=\frac18$, and  the above expression has value $(\frac{1}{16})(\frac{5}{32})$.

Now consider the case $\alpha_3\geq \frac18 +\frac{\alpha_2}2$. Then the right hand side of \eqref{eq:ex-bound} becomes
\[
\frac{1}{16}\left((\frac14+\alpha_2)(\frac18+\frac{\alpha_2}{2})-(\frac18+\frac{\alpha_2}{2})^2+ (\frac14+\alpha_3)\alpha_2-\alpha_2^2+\frac18-2\alpha_2\alpha_3\right).
\]
This expression is maximized (subject to our condition for this case)  when $\alpha_3=\frac18 +\frac{\alpha_2}2$
and $\alpha_2=\frac1{10}$. The value achieved is $(\frac{1}{16})(\frac{49}{320})$, and hence the previous case achieves the maximum. This gives an upper bound on $\Gamma(w,A)$ of $\frac{1}{16}\left(\frac{5}{32}\right)= \frac{5}{4}\left(\frac12\left(\frac14\right)^3\right)$. 
So we get $\frac{1}{2}(\frac14)^3\leq \Gamma(w)\leq \frac54(\frac12(\frac14)^3)$, as claimed. }
\end{remark}

\subsection{Black, white and grey regions in $\Delta$}\label{subse:region}
For the remainder of this section, fix $w\in\blue{W_0}$ with $\Gamma(w)>0$. Let  $\alpha=\Gamma(w)^{2/7}$, and let $R_w$ be the Robinson approximation of $w$ as described in Definition~\ref{def:R(w)}.
To prove Theorem~\ref{thm:main}, we need to break down $\Delta$  into regions where $R_w$ \blue{attains values between certain successive fractional multiples.}
We do so by having a closer look at the definition of $R_w$.

\begin{remark}\label{rem:regions}
Fix an integer $m$; we will determine the optimal value for $m$ later. 
For \blue{$1\leq k\leq m$} and $(x,y)\in \Delta$, the inequality $\frac{\blue{k-1}}{m}< R_w(x,y)\leq \frac{\blue{k}}{m}$ holds precisely when the following two conditions are satisfied:
\begin{itemize}
\item[(i)]  There exists an $\alpha\times \alpha$ cell $S\times T$ contained in $\UL(x,y)$ on which $\overline{w}(S,T)>\frac{\blue{k-1}}{m}$. 
\item[(ii)] For every $\alpha\times \alpha$ cell $S\times T$ contained in $\UL(x,y)$, we have  $\overline{w}(S,T)\leq \frac{\blue{k}}{m}$.
\end{itemize}

\end{remark}
Motivated by this observation, we use the graphon $w$ to split $\Delta$ into smaller regions of three types,  namely  black, white and grey regions as defined below. 
%
\begin{definition}\label{def:regions}
Let $\alpha=\Gamma(w)^{2/7}$, as was chosen in Definition~\ref{def:R(w)}. Let $m$ be an integer. 
For $k=\pink{1},\ldots, \blue{m}$, define the $k$'th black region $\B_k$, the $k$'th white region $\W_k$ and the $k$'th grey region $\G_k$ as follows. 
\begin{itemize}
\item $\B_k=\Big\{(x,y)\in \Delta:\ x=y\ \ \mbox{ or } \ \  \exists \ S\times T\subseteq \UL(x,y) \  \mbox{ with }\ |S|=|T|=\alpha \mbox{ and } \  \overline{w}(S, T)>\frac{\blue{k-1}}{m}\Big\}$. 
\item $\W_k=\Big\{(x,y)\in \Delta \setminus {\B}_k:\ \exists \ S\times T\subseteq \LR(x,y)\  \mbox{ with } \ |S|=|T|=\alpha\ \mbox{ and } \ \overline{w}(S, T)\leq \frac{\blue{k-1}}{m}\Big\}$.
\item $\G_k=\Delta\setminus (\B_k\cup\W_k)$.
\end{itemize}
Also, \pink{define $\B_0=\W_{m+1}=\Delta$ and $\W_0=\B_{m+1}=\emptyset$}.
\end{definition}

Finally, we introduce the regions of $\Delta$ on which the \blue{value of the Robinson approximation $R_w$ is easy to predict.}
\begin{definition}
For $0\leq k\leq m$, define $\R_k:=\B_k\cap \W_{k+1}$.
\end{definition}
See Figure~\ref{figure-regions} for a demonstration of the black, white and ${\mathcal R}_k$ regions. Note that no set $\B_k$, $k>0$, can contain any point within  \pink{distance} $\alpha$ of the border of  $[0,1]^2$; in the figure this margin is assumed to be invisible. A similar statement holds for $\W_k$ and points within  distance $\alpha$ of the diagonal, when $k\leq m$.
From Remark \ref{rem:regions}, it follows  that, \pink{for $1\leq k \leq m$,} $\B_{k}\setminus \B_{k+1}$ is exactly the region on which  $\frac{k-1}{m}<R_w\leq \frac{k}{m}$.
Note that $\R_k$ is a (possibly proper) subset of $\B_k\setminus \B_{k+1}$. \pink{This then leads to the following remark.
\begin{remark}
\label{rem:valuesR}
For $0\leq k\leq m$,
\[
\text{for all }(x,y)\in \R_k,\quad \,\frac{k-1}{m}<R_w(x,y)\leq \frac{k}{m}.
\]
For $k>0$ this follows immediately from the discussion above. For $k=0$, note that for all $1\leq i\leq m$,  $\R_0\cap (\B_i\setminus \B_{i+1})=\emptyset$ and thus $R_w$ cannot attain values in $\cup_{i=1}^m(\frac{i-1}{m},\frac{i}{m}]$ on $\R_0$. Since $R_w$ is bounded between 0 and 1, it follows that $R_w$ equals zero on $\R_0$. 
\end{remark}}

The intuition behind these definitions is the following.
If $(x,y)\in \R_k=\B_k\cap\W_{k+1}$, then $\UL(x,y)$  contains an $\alpha\times \alpha$ cell with ``high'' $w$-average, and $\LR(x,y)$ contains an $\alpha\times \alpha$ cell with ``low'' $w$-average. These two conditions, together with Lemma~\ref{lem:Gamma-lower-bound}, force us to assign a value to $R_w$
on the region $\R_k$, that lies between these two averages. For $(\B_k\setminus \B_{k+1})\setminus \R_k$, which lies in the grey region, we do not have such {\it a priori} knowledge, as Lemma~\ref{lem:Gamma-lower-bound} cannot be applied here. Luckily, this does not create problems for our error estimates, since the area of the grey regions turns out to be sufficiently small.

\begin{lemma}\label{lem:region-property}
With assumptions as in Definition~\ref{def:regions}, we have the following.
\begin{itemize}
\item[(i)] For every $k=0,\ldots, m$, we have $\B_{k+1}\subseteq \B_{k}$ and $\W_k\subseteq \W_{k+1}$. 
\item[(ii)] If $(x,y)\in \B_k$ then  $\LR(x,y)\subseteq \B_k$.  Similarly, if $(x,y)\in \W_k$ then $\UL(x,y)\subseteq \W_k$.
\item[(iii)] If $(x_1,y_1),(x_2,y_2)\in \G_k$ then $\UL(x_1,y_1)\cap \LR(x_2,y_2)\subseteq \G_k$.
\item[(iv)] For every $k=0,\ldots, m$, we have $\B_{k+1}\cap \W_{k}=\emptyset$.
\end{itemize}
\end{lemma}
\begin{proof}
Parts (i) \blue{and (ii) are easy to verify.}
For part (iv), note that 
$\B_{k+1}\cap \W_{k}\subseteq \B_k\cap \W_k$, and the latter set is empty by definition.

To prove (iii), suppose $(z,w) \in \UL(x_1,y_1)\cap \LR(x_2,y_2)\setminus \G_k$, which means that either 
$(z,w)\in \UL(x_1,y_1)\cap \LR(x_2,y_2)\cap \B_k$ or $(z,w)\in \UL(x_1,y_1)\cap \LR(x_2,y_2)\cap \W_k$. By (ii), the first case implies that $(x_1,y_1)\in \B_k$, and the second one implies that 
$(x_2,y_2)\in \W_k$. So either case leads to a contradiction. 
\end{proof}
\begin{lemma}\label{lem:properties-R}
For every $i\neq j$, we have $\R_i\cap \R_j=\emptyset$. Moreover, 
\begin{equation}\label{eq-R-i-partition}
\Delta\setminus \big(\bigcup_{k=1}^m \G_k\big)=\bigcup_{k=0}^{m} \R_k.
\end{equation}
\end{lemma}
\begin{proof}
By Definition~\ref{def:regions}, $\B_{k+1}\cap \W_{k+1}=\emptyset$, which implies that $\R_k\cap \R_{k+1}=\emptyset$.  Also from Lemma~\ref{lem:region-property} parts (i) and (iv), it follows that $\R_i\cap \R_j=\emptyset$ if $i\leq j-2$, since $\W_{i+1}\cap \B_j\subseteq \W_{i+1}\cap \B_{i+2}=\emptyset$. This proves that the regions $\R_i$ are disjoint. 
To show Equation (\ref{eq-R-i-partition}), observe that
\begin{eqnarray*}
\Delta\setminus \big(\bigcup_{k=1}^m \G_k\big)&=&\Delta\setminus \Big(\bigcup_{k=1}^m \ (\Delta\setminus(\W_k\cup \B_k))\Big)=\bigcap_{k=1}^m (\W_k\cup \B_k).
\end{eqnarray*}
Now, consider the expansion of $(\W_1\cup\B_1)\cap(\W_2\cup\B_2)\cap\ldots\cap(\W_m\cup\B_m)$ into expressions $X_1\cap \ldots \cap X_m$ with $X_i\in\{\W_i,\B_i\}$, and note that by Lemma~\ref{lem:region-property}, $X_1\cap \ldots \cap X_m=\emptyset$ whenever $X_i=\W_i$ and $X_j=\B_j$ for some $i<j$. So, every nonempty term $X_1\cap \ldots \cap X_m$  from the above expansion must be of one of the following forms:
\begin{itemize}
\item[(i)] $X_1\cap \ldots \cap X_m=\B_j\cap \W_{j+1}$ with $1\leq j<m$, if there is at least one black and one white region amongst $X_i$'s. 
\item[(ii)] $X_1\cap \ldots \cap X_m=\W_1\cap\cdots\cap\W_m=\W_1$, if all $X_i$'s are white. 
\item[(iii)] $X_1\cap \ldots \cap X_m=\B_1\cap\cdots\cap\B_m=\B_m$, if all $X_i$'s are black.
\end{itemize}
This finishes the proof, as $\B_m=\B_m\cap \W_{m+1}$ and $\W_1=\W_1\cap \B_0$.
\end{proof}
\begin{figure}
\centering
\begin{tikzpicture}[scale=0.33]
\draw (0,0) -- (10,0) -- (10,10) -- (0,10) -- (0,0);
\draw [fill=gray!10] (0,0)--(0,10)--(10,10)--(0,0);
\path [pattern= north east lines] (0,0)--(0,6) [out=70, in=190] to (4.2,10) --(10,10)--(0,0);
\path [pattern= north west lines] plot [smooth] coordinates {(0,5) (2,6) (4,7.1) (8,8.3) (9,9) (10,10)}
(10,10)--(0,0)--(0,5);
\draw [thick] (0,6) node[left,align=left] {\tiny $f_1\rightarrow$} (0,6) [out=70, in=190] to (4.2,10) ;
\draw [thick] (0,5) node[left,align=left] {\tiny $f_2\rightarrow$} plot [smooth] coordinates {(0,5) (2,6) (4,7.1) (8,8.3) (9,9) (10,10)};
\filldraw[fill=gray!10!white, draw=black] (3.5,-0.5) rectangle (4.5,-1.5); 
\node at (6.7,-1) {\small $\mathcal{B}_0=\Delta$};
\filldraw[pattern=north east lines, draw=black] (3.5,-1.7) rectangle (4.5,-2.7);
\node at (5.6,-2.2) {\small $\mathcal{B}_1$};
\filldraw[pattern= north west lines, draw=black] (3.5,-2.9) rectangle (4.5,-3.9);
\node at (5.6,-3.4) {\small $\mathcal{B}_2$};
\draw [thick, dashed]  (0,0) -- (10,10);
\draw[draw=white] (1.5,-4.1) rectangle (2.5,-5.4);
\end{tikzpicture}
\hspace{0.75cm}
\begin{tikzpicture}[scale=0.33]
\draw (0,0) -- (10,0) -- (10,10) -- (0,10) -- (0,0);
\draw [fill=gray!10] (0,0)--(0,10)--(10,10)--(0,0);
\path [pattern= horizontal lines] (0,10) -- (0,7.7) -- (1,8) [out=60, in=220] to (2,9.1) -- (2.7,10) --(0,10);
\path [pattern= vertical lines] plot [smooth] coordinates {(0,5.5) (3,7) (6,8) (8,9) (9,10)} (9,10)--(0,10)--(0,0)--(0,5.5);
\draw [thick] plot [smooth] coordinates {(0,5.5) (3,7) (6,8) (8,9) (9,10) }
(9,10) node [above] {\tiny $\begin{array}{c} g_2\\ \downarrow\end{array}$};
\draw [thick] (0,7.7) -- (1,8) [out=60, in=220] to (2,9.1) -- (2.7,10)
node [above] {\tiny $\begin{array}{c} g_1\\ \downarrow\end{array}$};
\filldraw[pattern=horizontal lines, draw=black] (3.5,-0.5) rectangle (4.5,-1.5); 
\node at (5.6,-1) {\small $\mathcal{W}_1$};
\filldraw[pattern=vertical lines, draw=black] (3.5,-1.7) rectangle (4.5,-2.7);
\node at (5.6,-2.2) {\small $\mathcal{W}_2$};
\filldraw[fill=gray!10, draw=black] (3.5,-2.9) rectangle (4.5,-3.9);
\node at (6.7,-3.4) {\small $\mathcal{W}_3=\Delta$};
\draw [thick, dashed]  (0,0) -- (10,10);
\draw[draw=white] (1.5,-4.1) rectangle (2.5,-5.4);
\end{tikzpicture}
%
\begin{tikzpicture}[scale=0.33]
\draw (0,0) -- (10,0) -- (10,10) -- (0,10) -- (0,0);
\fill[gray] plot [smooth] coordinates {(0,5.5) (3,7) (6,8) (8,9) (9,10) (10,10) (9,9) (8,8.3) (4,7.1) (2,6) (0,5)} (0,5)--(0,5.5);
\fill[gray] (0,6) [out=70, in=190] to (4.2,10) -- (2.7,10) -- (2,9.1) [in=220,out=60] to (1,8) -- (0,7.7) -- (0,6); 
\path [pattern= horizontal lines] (0,10) -- (0,7.7) -- (1,8) [out=60, in=220] to (2,9.1) -- (2.7,10) --(0,10);
\fill[gray!10] (0,10) -- (0,7.7) -- (1,8) [out=60, in=220] to (2,9.1) -- (2.7,10) --(0,10);
\path [pattern= horizontal lines] (0,10) -- (0,7.7) -- (1,8) [out=60, in=220] to (2,9.1) -- (2.7,10) --(0,10);
\path [pattern= north east lines] (0,5.5)--(0,6) [out=70, in=190] to (4.2,10) --(--(9,10) plot [smooth] coordinates {(9,10) (8,9) (6,8) (3,7) (0,5.5)};
\path [pattern= vertical lines] (0,5.5)--(0,6) [out=70, in=190] to (4.2,10) --(--(9,10) plot [smooth] coordinates {(9,10) (8,9) (6,8) (3,7) (0,5.5)};
\fill[gray!10]plot [smooth] coordinates {(0,5) (2,6) (4,7.1) (8,8.3) (9,9) (10,10)}
(10,10)--(0,0)--(0,5);
\path [pattern= north west lines] plot [smooth] coordinates {(0,5) (2,6) (4,7.1) (8,8.3) (9,9) (10,10)}
(10,10)--(0,0)--(0,5);
\draw [thick] (0,7.7) -- (1,8) [out=60, in=220] to (2,9.1) -- (2.7,10) node [above] {\tiny $\begin{array}{c} g_1\\ \downarrow\end{array}$};
\draw [thick] (0,6) node[left,align=left] {\tiny $f_1\rightarrow$} (0,6) [out=70, in=190] to (4.2,10) ;
\draw [thick] plot [smooth] coordinates {(0,5.5) (3,7) (6,8) (8,9) (9,10) }
(9,10) node [above] {\tiny $\begin{array}{c} g_2\\ \downarrow\end{array}$};
\draw [thick] (0,5) node[left,align=left] {\tiny $f_2\rightarrow$} plot [smooth] coordinates {(0,5) (2,6) (4,7.1) (8,8.3) (9,9) (10,10)};
\draw [thick] plot [smooth] coordinates {(0,5) (2,6) (4,7.1) (8,8.3) (9,9) (10,10)};
\fill[gray!10] (1.5,-0.5) rectangle (2.5,-1.5); 
\filldraw[pattern=horizontal lines, draw=black] (1.5,-0.5) rectangle (2.5,-1.5);
\node at (6,-1) {\small $\mathcal{R}_0=\mathcal{B}_0\cap \mathcal{W}_1$};
\fill[pattern=north east lines] (1.5,-1.7) rectangle (2.5,-2.7);
\filldraw[pattern=vertical lines, draw=black] (1.5,-1.7) rectangle (2.5,-2.7);
\node at (6,-2.2) {\small $\mathcal{R}_1=\mathcal{B}_1\cap\mathcal{W}_2$};
\fill[fill=gray!10] (1.5,-2.9) rectangle (2.5,-3.9);
\filldraw[pattern= north west lines, draw=black] (1.5,-2.9) rectangle (2.5,-3.9);
\node at (6,-3.4) {\small $\mathcal{R}_2=\mathcal{B}_2\cap\mathcal{W}_3$};
\draw [thick, dashed]  (0,0) -- (10,10);
\filldraw[fill=gray, draw=black] (1.5,-4.1) rectangle (2.5,-5.1);
\node at (6,-4.7) {\small Grey regions};
\draw [thick, dashed]  (0,0) -- (10,10);
\end{tikzpicture}
\caption{Black and white regions; example for $m=2$.}
\label{figure-regions}
\end{figure}
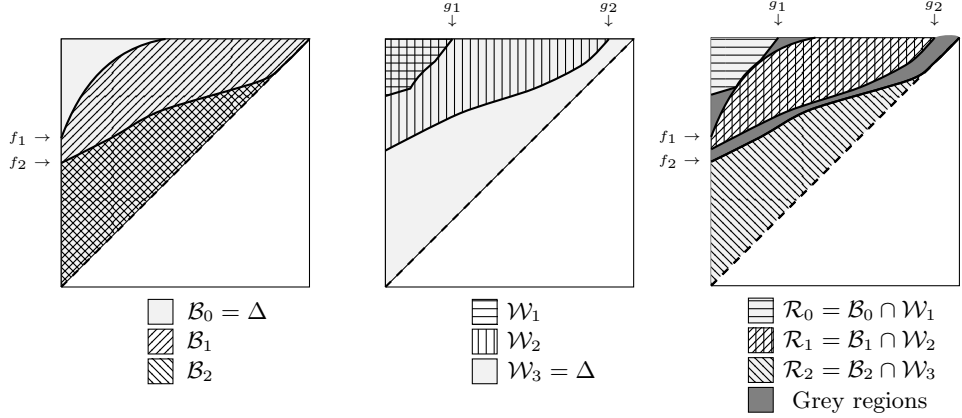
\begin{remark}\label{rem:upper-lower-boundary}
For every $1\leq k\leq m$, the region $\G_k$  is bounded between \emph{lower} and   \emph{upper boundary} functions $f_k,g_k:[0,1]\rightarrow [0,1]$, as shown in Figure~\ref{figure-regions}. Indeed, $f_k$ is the upper  boundary of $\B_k$, and $g_k$ is the lower boundary of $\W_{k}$. To be precise, for every $x\in[0,1]$, 
\begin{eqnarray}
f_k(x)&=&\sup\{z\in[x,1]: (x,z)\in \B_k\},\label{eq:f-function}\\
g_k(x)&=&\inf\{z\in[x,1]: (x,z)\in \W_k\}, \label{eq:g-function}
\end{eqnarray}
where we set $\inf\emptyset =1$, if it appears in the definition of $g_k$. In addition, we define $f_0(x)= 1$ and $g_{m+1}(x)=x$ for all $x\in[0,1]$, to represent the corresponding boundaries for $\B_0=\W_{m+1}=\Delta$. \blue{Finally, since $f_k$ and $g_{k+1}$ are the upper and lower boundaries of $\B_k$ and $\W_{k+1}$ respectively, if the region  $\R_k$ is nonempty, then it is  bounded from below by  $g_{k+1}$, and from above by $f_{k}$. }

From the definition of $\B_k$ and $\W_k$, it is clear that $f_k\leq g_k$.
Moreover, both $f_k$ and $g_k$ are increasing functions. Indeed, towards a contradiction suppose $f_k$ is not increasing,  i.e.~assume that there exist $x_1,x_2\in [0,1]$, with $x_1<x_2$ but $f_k(x_1)>f_k(x_2)$. Then, 
$(x_2,\frac{f_k(x_1)+f_k(x_2)}{2})\not\in \B_k$ but $(x_2,\frac{f_k(x_1)+f_k(x_2)}{2})\in \LR(x_1,\frac{f_k(x_1)+f_k(x_2)}{2})$, which is a contradiction with  $(x_1,\frac{f_k(x_1)+f_k(x_2)}{2})\in \B_k$ by Lemma~\ref{lem:region-property} (ii). The proof for $g_k$ is similar. 

Note that boundary functions $f_k, g_k$ are not necessarily continuous. However, being increasing, the boundary functions can only admit jump discontinuities. So, we can naturally extend the graph of a boundary function to a curve by adding appropriate vertical line segments at its points of discontinuity. We call the resulting curve a \emph{boundary curve}, and we denote it by $f_k,g_k$ again. 
\end{remark}
%

\begin{lemma}\label{lem:grey-region-small}
Assume the notations and conditions of Definition~\ref{def:regions}, and let $k\in\{1,\ldots, m\}$. Then,
$\overline{\G_k}$ does not contain any $\beta\times \beta$ square, where $\beta>\alpha$. Here, $\overline{\G_k}$ denotes the closure of $\G_k$ in ${\mathbb R}^2$. 
\end{lemma}
\begin{proof}
Let $k\in\{1,\ldots, m\}$ be fixed. Clearly, points of $\overline{\G_k}$ lie on or between the lower and upper boundary curves  of $\G_k$. 
Towards a contradiction, suppose there are $\beta>\alpha$, and measurable subsets $S,T\subseteq [0,1]$ with  $|S|=|T|=\beta$ for which $S\times T\subseteq \overline{\G_k}$. Let $a_1=\inf S$, $a_2=\sup S$, $b_1=\inf T$, and $b_2=\sup T$; note that $a_2-a_1\geq \beta$ and $b_2-b_1\geq \beta$. 
Since $\overline{\G_k}$ is closed, we have $(a_i,b_j)\in \overline{\G_k}$ for $i,j=1,2$; thus  using Lemma~\ref{lem:region-property}~(iii), one can easily see that $[a_1,a_2]\times [b_1,b_2]\subseteq \overline{\G_k}$. Note that every point in $\overline{\G_k}$, that is not on the lower or upper boundary curves, must be an inner point. Given that the lower and upper boundary functions of $\G_k$ are increasing, we conclude that $(a_1,a_2)\times (b_1,b_2)\subseteq \G_k$. Clearly, $(a_1,a_2)\times (b_1,b_2)$ contains a closed $\alpha\times \alpha$ rectangle, say 
$[a_1',a_2']\times [b_1',b_2']$. 

The two points $(a'_1,b'_2), (a'_2,b'_1)$ belong to $\G_k$, so both of them fail to satisfy the conditions for $\W_k$ and $\B_k$. In particular, we have
$\overline{w}([a_1',a_2'], [b_1',b_2'])\leq \frac{\blue{k-1}}{m}$ as $(a'_2,b'_1)\not\in \B_k$, and $\overline{w}([a_1',a_2'], [b_1',b_2'])>\frac{\blue{k-1}}{m}$ as $(a'_1,b'_2)\not\in \W_k$.  This is a contradiction. 
\end{proof}

\subsection{Proof of Theorem~\ref{thm:main}.}
We are now ready to prove Theorem~\ref{thm:main}. We use the contrapositive, by showing that, if $\| R_w-w\|_{\Box}$ is large, then $\Gamma (w)$ cannot be small.  
First observe that if $\Gamma (w)^{1/7}>1/14$, then $\| R_w-w\|_{\Box}\leq 1<14\Gamma(w)^{1/7}$, and there is nothing to prove. So without loss of generality assume that $0<\Gamma (w)^{1/7}\leq 1/14$. 

Let $m$ be an integer chosen from $[\frac{13}{14}\Gamma(w)^{-1/7}, \Gamma(w)^{-1/7}]$. Note that such an integer exists, as the interval has length at least 1. Moreover, from the upper bound on $\Gamma (w)$, we get $m\geq 13$.   Let $\alpha$ be as in Definition~\ref{def:R(w)},  i.e.  $\alpha= \Gamma(w)^{2/7}$. 
Note that we chose these parameters so that $\frac{\alpha^3}{m}\geq \Gamma(w)$. This choice of parameters allows us to use  Lemma~\ref{lem:Gamma-lower-bound} to obtain a contradiction. In the rest of the proof, we will show that  if $\|w-R_w\|_\Box$ is not small, then we can find two rectangles satisfying  conditions of Lemma~\ref{lem:Gamma-lower-bound}, which will result in a contradiction. 

Fix $\beta>\alpha$, and let $\delta:=14 m\beta$. 
Towards a contradiction, assume that $\|w-R_w\|_\Box> \delta$. So there exist measurable subsets $S,T\subseteq [0,1]$ so that
\begin{equation*}\label{eq;choice-D}
\Big|\int\int_{S\times T} w-R_w\Big|>\delta.
\end{equation*} 
Replacing $ S\times T$ with $T\times S$ if necessary, we can assume \emph{wlog} that $\Big|\int\int_{(S\times T)\cap \Delta} w-R_w\Big| \, >\frac{\delta}{2}$. \pink{We will now show that $(S\times T)\cap\Delta$ contains a region $S^*\times T^*$ so that $w$ is substantially different from $R_w$ on this region. The precise statement is given in the following claim.
\begin{claim}\label{claim:S*T*}
There exist sets $S^*,T^*$ of size $|S^*|=|T^*|=\alpha$ and $0\leq k\leq m$
so that $S^*\times T^*\subseteq \R_k$, and for which we have the following inequality:
\begin{equation}\label{eq-proof-fix}
\Big| \int\int_{S^*\times T^*} w-R_w\Big|\, \geq \,  m\alpha^3 (3-\frac{4}{m}).
\end{equation}
\end{claim}}
\begin{proof}
Let $S$ and $T$ be as described above. 
Split $S$ into $N_1=\lceil |S|/\beta\rceil$ subsets $S_1\leq S_2 \leq\dots \leq  S_{N_1}$ so that $|S_1|=|S_2|=\dots=|S_{N_1-1}|=\beta$ and $|S_{N_1}|\leq \beta$.
Likewise, we split $T$ into $N_2=\lceil |T|/\beta\rceil$ subsets $T_1\leq T_2\leq \dots \leq T_{N_2}$ with $|T_1|=|T_2|=\dots=|T_{N_2-1}|=\beta$ and $|T_{N_2}|\leq \beta$. 
\blue{The sets $S_i,T_j$ form a ``grid'', in which each cell is of the form $S_i\times T_j$ with $1\leq i\leq N_1$ and $1\leq j\leq N_2$.}
Since  $|S_{N_1}|\leq \beta$, $|T_{N_2}|\leq \beta$, and $|w-R_w|\leq 1$, we can use triangle inequality to get
\begin{eqnarray}
\sum_{\scriptsize{
\begin{array}{c}
1\leq i< N_1, \ 1\leq j< N_2\\
(S_i\times T_j)\cap \Delta\neq \emptyset\\
\end{array}}}
\left|\int\int_{S_i\times T_j} w-R_w\right|
&\geq&
\left|\sum_{\scriptsize{
\begin{array}{c}
1\leq i< N_1, \ 1\leq j< N_2\\
(S_i\times T_j)\cap \Delta\neq \emptyset\\
\end{array}}}
\int\int_{S_i\times T_j} w-R_w\right|\nonumber\\
&\geq& \left|\int\int_{(S\times T)\cap \Delta} w-R_w\right|-\left|\int\int_{(S_{N_1}\times T)\cup (S\times T_{N_2})} w-R_w\right|\nonumber\\
&>&\frac{\delta}{2}-2\beta. \label{eq:int-over-sq-grid}
\end{eqnarray} 
Recall that $\Delta=\bigcup_{k=0}^m\R_k\, \cup  \, \bigcup_{k=1}^m \G_k$, and each of the regions $\G_k$ or $\R_k$ is bounded by boundary curves  from the collection 
$\{f_k, g_l: \ 1\leq k\leq m, 1\leq l\leq m+1\}$ as defined in Remark~\ref{rem:upper-lower-boundary}. Thus, if a  cell $S_i\times T_j$  does not cross the graph of any of these boundary curves, then it must be entirely contained inside one closed region $\overline{\R_k}$ or $\overline{\G_k}$. Here,  by ``a cell crossing a boundary'', we mean that the  top-left corner of the cell is strictly above the boundary curve, and its bottom-right corner is  strictly below the curve. (Note that  we need to use the concept of ``a cell crossing a boundary'' rather than ``a cell intersecting with the boundary'', as our cells are not necessarily connected subsets of ${\mathbb R}^2$. So a boundary curve can go through a cell, without having to intersect with it.)
Next,  by Lemma~\ref{lem:grey-region-small}, none of the  grey regions $\overline{{\G_k}}$ can contain any of the cells ${S_i}\times {T_j}$ with $1\leq i<N_1$ and $1\leq j<N_2$.
Let $\I$ denote the collection of indices $(i,j)$ with $i<N_1$ and $j<N_2$, for  which the associated cells do not lie in a single region $\overline{\R_k}$. From the above discussion, we have
$$\I=\Big\{(i,j): 1\leq i<N_1, 1\leq j< N_2,\ \mbox{ and } \exists\  1\leq k \leq m \mbox{ s.t. } (S_i\times  T_j) \mbox{ crosses } f_k \mbox{ or } g_{k} \mbox{ or }  g_{m+1}\Big\}.$$
\pink{Note that any cell with indices outside $\I$ is completely contained in the closure of one region $\R_k$.}

Now, we  bound the size of the set $\I$.
By Remark~\ref{rem:upper-lower-boundary}, the lower and upper boundaries $f_k,g_k$  are increasing functions. We claim that $f_k$ (similarly $g_k$) crosses at most $2/\beta$ cells from the grid. Indeed, suppose that  $S_{n_1}\times T_{n'_1}, S_{n_2}\times T_{n'_2},\ldots, S_{n_{p}}\times T_{n'_{p}}$ is a sequence of distinct cells, all of which cross the graph of $f_k$. Since $f_k$ is increasing, after relabeling if necessary, we have that $n_1\leq n_2\leq \ldots\leq n_p$ and $n'_1\leq n'_2\leq\ldots\leq n'_p$. However, $\Delta$ cannot contain any such sequence of length more than $N_1-1+N_2-1\leq 2/\beta$. Thus, we have 
$$|\I|\leq \frac{2(2m+1)}{\beta}.$$
Since every cell indexed in $\I$ is of size $\beta^2$, we have $|\F|\leq (4m+2)\beta$, where $\F=\bigcup_{(i,j)\in \I} S_i\times T_j$.
By inequality (\ref{eq:int-over-sq-grid}), and the fact that $|w-R_w|\leq 1$, we get
\begin{equation*}
\sum_{\scriptsize{
\begin{array}{c}
1\leq i< N_1, \ 1\leq j< N_2\\
(S_i\times T_j)\cap \Delta\neq \emptyset\\
(i,j)\not\in {\cal I}
\end{array}}}
\left|\int\int_{S_i\times T_j} w-R_w\right|>\frac{\delta}{2}-2\beta- (4m+2)\beta=3m\beta-4\beta.
\end{equation*}
By the pigeonhole principle, since there are at most $(1/\beta)^2$  cells  $S_i\times T_j$ of size $\beta\times \beta$, there must exist a cell $S_{i_0}\times T_{j_0}\subseteq \Delta\setminus \F$ so that $|S_{i_0}|=|T_{j_0}|=\beta$ and
\begin{equation*}
\Big| \int\int_{S_{i_0}\times T_{j_0}} w-R_w\Big| \geq m\beta^3(3-\frac{4}{m}).
\end{equation*}
\pink{Since $(i_0,j_0)\not\in \I$, we conclude that} $S_{i_0}\times T_{j_0}$ lies entirely \pink{in the closure}  $\overline{\R_k}=\overline{\B_k\cap \W_{k+1}}$ for some $0\leq k\leq m$. \blue{It is easy to see that, there are subsets  $S^*\subseteq S_{i_0}$ and $T^*\subseteq T_{j_0}$ with $|S^*|=|T^*|=\alpha$,
%
}
\textcolor{black}{so that} $S^*\times T^*\subseteq \R_k=\B_k\cap \W_{k+1}$ and
\begin{equation}\label{eq-proof-fix}
\Big| \int\int_{S^*\times T^*} w-R_w\Big| \geq m\alpha^2\beta(3-\frac{4}{m})\geq m\alpha^3 (3-\frac{4}{m}),
\end{equation}
which proves the claim.
\textcolor{black}{(In the case where $\frac{\alpha}{\beta}$ is a rational number, the existence of such sets $S^*$ and $T^*$ is an easy application of the pigeonhole principle and the fact that 
 $\frac{1}{\beta^2}\Big| \int\int_{S_{i_0}\times T_{j_0}} w-R_w\Big| \geq m\beta(3-\frac{4}{m}).$ 
 The case where $\frac{\alpha}{\beta}$ is irrational follows from a standard density argument.)}
\end{proof}

\pink{We now state and prove a second claim.
\begin{claim}
If sets $S^*,T^*$ and index $k$ as in Claim \ref{claim:S*T*} exist, then $\Gamma (w)>\alpha^3/m$.
\end{claim}}
\begin{proof} 
\pink{As explained in Remark \ref{rem:valuesR}}, we have $\frac{k-1}{m}< R_w\leq\frac{k}{m}$ on \blue{$S^*\times T^*$}. 
\pink{To prove the claim, we will consider two cases.
 }

\pink{{\bf Case 1.}} Assume first that $w$ has larger average than $R_w$ on \blue{$S^*\times T^*$}, so $\int\int_{\blue{S^*\times T^*}} w-R_w\geq m\alpha^3(3-\frac{4}{m})$. 
In this case, we have
\begin{equation}\label{eqn:avg}
\overline{w}(\blue{S^*, T^*})-\frac{\blue{k-1}}{m}\blue{>} \overline{w-R_w}(\blue{S^*, T^*})\geq \frac{m\alpha^3(3-\frac{4}{m})}{|S^*\times T^*|}=m\alpha(3-\frac{4}{m}),
\end{equation}
where we use that $|S^*\times T^*|=\alpha^2$.


\pink{Next we argue that $k<m$. From our choice of parameters, we have that $m^2\alpha\geq (\frac{13}{14})^2>\frac{1}{2}$, which implies that $2m\alpha-\frac{1}{m}>0$. Assuming now that $S^*\times T^*\subseteq \R_m$, 
then by \eqref{eqn:avg} and the facts that $0\leq w\leq 1$ and $m\geq 13$, we have that
$$
1<1+(2m\alpha-\frac{1}{m})<(1-\frac{1}{m})+m\alpha(3-\frac{4}{m})<\overline{w}(S^*, T^*)\leq 1,
$$
which is a contradiction. This proves that $k<m$.}
 
Now let $(x,y)$ be the lower right corner of $S^*\times T^*$, so $x=\sup S^*$ and $y=\inf T^*$. Then $(x,y)\in \W_{k+1}$, and thus $\LR (x,y)$ contains a region $S_l\times T_l$ so that $|S_l|=|T_l|=\alpha$, and $\overline{w}(S_l, T_l)\leq  \blue{\frac{k}{m}}$. Applying Lemma~\ref{lem:Gamma-lower-bound} together with  inequality (\ref{eqn:avg}), we now conclude that 
$$\Gamma(w)\geq \alpha^3\big(m\alpha(3-\frac{4}{m})+ \frac{k-1}{m}-\frac{k}{m}\big).$$
Since  $m\geq\max\{ \frac{13}{14} \Gamma(w)^{-1/7}, 13\}$, we conclude that
\begin{equation}\label{eq:bounds}
\Gamma(w)(\frac{m}{\alpha^3})\geq  m^2\alpha(3-\frac{4}{m})-1\geq (\frac{13}{14})^2\Gamma(w)^{-2/7}\Gamma(w)^{2/7}(3-\frac{4}{13})-1>1,
\end{equation}
Therefore,
\[
\Gamma (w)> \frac{\alpha^3}{m},
\]
\pink{which proves the claim for the first case.}

\pink{{\bf Case 2.}} For the second case, assume $\int\int_{\blue{S^*\times T^*}}R_w-w\geq m\blue{\alpha}^3(3-\frac{4}{m})$. 
\pink{This case can only happen if $k>0$, as $R_w$ attains 0 on $S^*\times T^*$ if $k=0$.}
By a similar argument,
\begin{equation}\label{eq:bounds2}
\frac{\blue{k}}{m}-\overline{w}(\blue{S^*, T^*})\geq m\blue{\alpha}(3-\frac{4}{m}).
\end{equation}

Now let $(x,y)$ be the upper left corner of $\blue{S^*\times T^*}$, so $x=\inf \blue{S^*}$ and $y=\sup \blue{T^*}$. Then $(x,y)\in \B_{k}$, and thus $\UL (x,y)$ contains a region $S_l\times T_l$ so that $|S_l|=|T_l|=\alpha$, and $\overline{w}(S_l, T_l)> \frac{\blue{k-1}}{m}$. Applying Lemma~\ref{lem:Gamma-lower-bound} together with  \blue{Equations (\ref{eq:bounds2}) and \eqref{eq:bounds}}, we now conclude that 
$$\Gamma(w)\geq \alpha^3\big(m\blue{\alpha}(3-\frac{4}{m})-\frac{1}{m}\big)=\frac{\alpha^3}{m}\big(m^2\blue{\alpha}(3-\frac{4}{m})-1\big)  >\pink{\frac{\alpha^3}{m}}.
$$
\end{proof}

\pink{Our second claim directly contradicts the assumption that  $\Gamma(w)\leq \frac{\alpha^3}{m}$. Thus our assumption that  $\| R_w-w\|_{\Box}>\delta $ does not hold, which implies that $\| R_w-w\|_{\Box}\leq \delta =14m\beta$ for all $\beta >\alpha$. Thus, $\| R_w-w\|_{\Box}\leq 14m\alpha\leq 14\Gamma(w)^{1/7}$. This completes the proof of Theorem \ref{thm:main}.}
\section{Recognition of graph sequences sampled from Robinson graphons}\label{Sec:conjecture}
Consider a  graph sequence $\{G_n\}$ sampled from a graphon $w\in \blue{W_0}$, i.e. a sequence $\{G_n\}$ that converges to the graphon $w$ in the sense of Lov\'{a}sz-Szegedy, or equivalently $\delta_\Box(w_{G_n},w)\rightarrow 0$. 
By \cite[Theorem 6.4]{CGHJK}, $\widetilde{\Gamma}$ is continuous, and in particular,
\begin{equation}\label{easy-direction}
\mbox{ If } \ G_n\rightarrow w\ \mbox{ and } \ \Gamma(w)=0, \ \mbox{ then } \ \widetilde{\Gamma}(G_n)\rightarrow 0. \mbox{\tag{$\star$}}
\end{equation} 
Continuity guarantees a weaker version of the converse: if $G_n\rightarrow w$ and $\widetilde{\Gamma}(G_n)\rightarrow 0$, then 
$\widetilde{\Gamma}(w)=0$. Since $\widetilde{\Gamma}(w)=\inf\{ \Gamma (u):\delta_\Box (w,u)\}$, these earlier results do not imply that the $\delta_\Box$-equivalence class of $w$ contains a Robinson graphon. In Theorem~\ref{thm-conj} to follow, we will use Theorem~\ref{thm:main} and the weak$^*$  topology of the space of graphons to show the strong version of the converse.
The proof of Theorem~\ref{thm-conj} involves approximations of a graphon by step graphons, which we discuss in the following subsection.  

\subsection{Stepping operator}
For an integer $N\in\bbN$ and a graphon $u\in \blue{W_0}$, we define the step graphon $u^{(N)}$ as follows: split the interval $[0,1]$ into $N$ equal-sized subintervals $I_1, \dots, I_{N}$, and define
$$u^{(N)}(x,y)=\blue{\overline{u}({I_i,  I_j}),} \mbox{ if } (x,y)\in I_i\times I_j.$$
To avoid clutter of notation, \blue{we use $u_{i,j}^{(N)}$  to denote 
the value of $u^{(N)}$ at every point in $I_i\times I_j$.} 
The operator that assigns to every $u\in \blue{W_0}$, the step graphon $u^{(N)}$ is called a \emph{stepping operator.}
Step graphons approximate $u$ in $\|\cdot\|_1$-norm, and hence in $\|\cdot\|_\Box$-norm. In the following lemma, we obtain a universal upper bound for the rate of convergence of the step graphon approximation of a Robinson graphon.
Note that the proof of Lemma~\ref{lem:stepping} relies heavily on the Robinson structure of the graphon. For a general graphon, a uniform bound for the rate of convergence of step graphon approximations can be obtained by applying the (Weak) Regularity Lemma for graphons (for example, see \cite[Lemma 9.9]{lovasz-book}). However, the bound from the Regularity Lemma cannot be used for our methods of proving Theorem~\ref{thm-conj}, since the bound is of the order of $\frac{1}{\sqrt{N}}$, and the partition sets are not necessarily intervals, but rather measurable sets. The following lemma gives us better control on both the partition sets and the error bound. 
\begin{lemma}\label{lem:stepping}
Let $u\in \blue{W_0}$ be a Robinson graphon. Then for every $N\in {\mathbb N}$, $\|u-u^{(N)}\|_1\leq \frac{7}{N}$
\end{lemma}
\begin{proof}
Fix $N\in \bbN$, and note that $u^{(N)}$ is a Robinson graphon as well. Consider the symmetric graphons $u_-^{(N)}$ and $u_+^{(N)}$  that are obtained from $u^{(N)}$ by shifting every cell towards the diagonal (for $u_-^{(N)}$) or away from the diagonal (for $u_+^{(N)}$).
That is, 
\begin{eqnarray*}
u_-^{(N)}(x,y):=u_-^{(N)} (y,x)=
\left\{
\begin{array}{cl}
u^{(N)}_{i-1,j+1}& \mbox{ if } (x,y)\in I_i\times I_j \ \mbox{ and }\  1<i\leq j<N  \\
0  & \mbox{ if }  (x,y)\in I_{i}\times I_{N}\   \mbox{ or }\  (x,y)\in I_1\times I_{j}
\end{array}
\right.
\end{eqnarray*}
and 
\begin{eqnarray*}
u_+^{(N)}(x,y):=u_+^{(N)} (y,x)=
\left\{
\begin{array}{cl}
u^{(N)}_{i+1,j-1}& \mbox{ if } (x,y)\in I_i\times I_j,\  1\leq i<j\leq N \ \mbox{ and } j-i\geq 2\\
1 & \mbox{ if }  (x,y)\in I_{i}\times I_{j} \  \mbox{ and }\  0\leq j-i\leq 1\\
\end{array}
\right.
\end{eqnarray*}
Since $u$ is a Robinson graphon, it can be easily checked that $u_-^{(N)}\leq u\leq u^{(N)}_+$ and $u_-^{(N)}\leq u^{(N)}\leq u^{(N)}_+$. On the other hand, 
\begin{eqnarray*}
\|u_+^{(N)}-u_-^{(N)}\|_1&=&\int_{[0,1]^2} u_+^{(N)}-u_-^{(N)}\\
&=&\frac{1}{N^2}\big(2\sum_{i=1}^{N-2}\sum_{j=i+2}^N u_{i+1,j-1}^{(N)}+\sum_{i=1}^N\sum_{j=i-1}^{i+1}1-2\sum_{i=2}^{N-1}\sum_{j=i+1}^{N-1} u_{i-1,j+1}^{(N)}-\sum_{i=2}^{N-1} u_{i-1,i+1}^{(N)}\big)\\
&\leq&\frac{7N}{N^2},
\end{eqnarray*}
because every pair $(i,j)$ with $1<i<j<N$ and $j-i\geq 3$, contributes exactly once in a positive sum and exactly once in a negative sum, and there are only $7N$ other terms left. Thus, we have 
$\|u-u^{(N)}\|_1\leq \|u_+^{(N)}-u_-^{(N)}\|_1\leq \frac{7}{N}$.
\end{proof}
\subsection{Robinsonian graphons}

We are now ready to  prove Theorem~\ref{thm-app}, which will follow from Theorem~\ref{thm-conj}. First  we introduce the following definition, which matches a similar concept in matrix theory. 
\begin{definition}
A graphon $w\in \blue{W_0}$ is called \emph{Robinsonian} if there exists a Robinson graphon $u\in \blue{W_0}$ such that $\delta_\Box(u,w)=0$. In other words, a graphon is Robinsonian if its $\delta_\Box$-equivalence class contains a Robinson graphon.  
\end{definition}
%

\begin{theorem}\label{thm-conj}
Let $\{G_n\}_{n\in{\mathbb N}}$ be a growing sequence of graphs converging to a graphon $w\in \blue{W_0}$. Then,
$w$ is Robinsonian if and only if $\widetilde{\Gamma}(G_n)\rightarrow 0$. 
\end{theorem}
\begin{proof}
The forward direction is a consequence of continuity of $\widetilde{\Gamma}$. 
To prove the backward direction, suppose $\widetilde{\Gamma}(G_n)\rightarrow 0$. Without loss of generality, assume that every $G_n$ is labeled so that $\widetilde{\Gamma}(G_n)$ is achieved, that is $\widetilde{\Gamma}(G_n)=\Gamma(w_{G_n})$ where $w_{G_n}$ denotes the graphon that represents $G_n$. From the assumption we have $\Gamma(w_{G_n})\rightarrow 0$, and from the definition of convergence of graph sequences, we have that $\delta_\Box (w_{G_n},w)\rightarrow 0$.

Applying Theorem~\ref{thm:main}, for every $n\in {\mathbb N}$, there exists a Robinson graphon $u_n\in \blue{W_0}$ such that  $\|u_n-w_{G_n}\|_\Box\leq 14\Gamma(w_{G_n})^{1/7}$. So $\{u_n\}_{n\in{\mathbb N}}$ is a sequence of Robinson graphons such that $\delta_\Box(u_n,w)\rightarrow 0$ as $n\rightarrow \infty$, because
$$\delta_\Box(u_n,w)\leq \delta_\Box(u_n,w_{G_n})+\delta_\Box(w_{G_n},w) \leq \|u_n-w_{G_n}\|_\Box+\delta_\Box(w_{G_n},w).$$

Now consider the graphon space $\blue{W_0}$ as a subset of $B_1(L^\infty[0,1]^2):=\left\{f\in L^\infty[0,1]^2: \ \|f\|_\infty\leq 1\right\}$, namely the closed unit ball of $L^\infty[0,1]^2$. The Banach space $L^\infty[0,1]^2$ is isometrically isomorphic to the (Banach space) dual of $L^1[0,1]^2$, so one can equip $B_1(L^\infty[0,1]^2)$ with the weak* topology induced by this duality. By the Banach-Alaoglu theorem, $B_1(L^\infty[0,1]^2)$ is compact in the weak* topology 
(\cite[Theorem 3.1 of Chapter V]{conway-fnl}). In addition, since $L^1[0,1]^2$ is separable, the unit ball  $B_1(L^\infty[0,1]^2)$ is a metrizable space in the weak* topology, and thus it is sequentially compact as well 
(\cite[Theorem 5.1 of V]{conway-fnl}). 
Thus the sequence $\{u_n\}_{n\in{\mathbb N}}$ in  $B_1(L^\infty[0,1]^2)$ has a weak* convergent subsequence. 

By going down to that subsequence if necessary, \emph{wlog} we can assume that 
$\{u_n\}_{n\in \bbN}$ converges to some $z\in B_1(L^\infty[0,1]^2)$ in the weak* topology, i.e.~for every $h\in L^1[0,1]^2$ we have  $\int_{[0,1]^2} u_nh\rightarrow \int_{[0,1]^2}zh$ as $n\rightarrow \infty$.  
In particular, for every measurable subsets $S,T\subseteq [0,1]$, we have
\begin{eqnarray}\label{weak*}
\int_{S\times T}u_n\rightarrow \int_{S\times T} z, \mbox{ as }n\rightarrow \infty.
\end{eqnarray}
By (\ref{weak*}), we get that  for every $N\in{\mathbb N}$, the sequence $\{u_n^{(N)}\}$ converges point-wise to $z^{(N)}$. Note that every graphon $u_n$ is Robinson. So for every $N\in{\mathbb N}$,  the corresponding step graphons $u_n^{(N)}$ and their point-wise limit $z^{(N)}$ are also Robinson. Finally, since $\|z- z^{(N)}\|_1\rightarrow 0$,  $z$ is a Robinson graphon a.e. as well. 

Next, we claim that $\delta_\Box(z,w)=0$. Fix $\epsilon>0$, and choose $N_0\in\bbN$ such that $\frac{7}{N_0}\leq \frac{\epsilon}{3}$. From the convergence of $\{u_n^{(N_0)}\}_{n\in \bbN}$ to $z^{(N_0)}$, pick $m\in \bbN$ so that $\|u_m^{(N_0)}-z^{(N_0)}\|_\infty\leq \frac{\epsilon}{3}$. 
Applying Lemma~\ref{lem:stepping} to the Robinson graphons $z^{(N_0)}$ and $u_m^{(N_0)}$, we have
\begin{equation}\label{eq-1}
\|u_m-z\|_1\leq \|u_m-u_m^{(N_0)}\|_1+\|u_m^{(N_0)}-z^{(N_0)}\|_1+\|z^{(N_0)}-z\|_1\leq \frac{7}{N_0}+\frac{\epsilon}{3}+\frac{7}{N_0}\leq \epsilon.
\end{equation}
Since $\delta_\Box(w,z)\leq \delta_\Box(w,u_m)+\delta_\Box(u_m,z)\leq \delta_\Box(w,u_m)+\|u_m-z\|_1$ for all $m$, and $\delta_\Box(w,u_m)+\|u_m-z\|_1\rightarrow 0$, we have that $\delta_\Box (w,z)=0$.
\end{proof}
\begin{corollary}
$\widetilde\Gamma(w) = 0$ if and only if $w$ is Robinsonian.
\end{corollary}
\begin{proof}
Note that $\widetilde{\Gamma}(w)$ is defined as an infimum over all graphons in the $\delta_\Box$-equivalence class of $w$. The above theorem tells us that this infimum is in fact achieved at 0. 
\end{proof}
\begin{remark}
\begin{itemize}
\item[(i)] A $\delta_\Box$-equivalence class of graphons may include more than one Robinson graphon, i.e.~there is no concept of a ``unique Robinson representation'' of a $\delta_\Box$-equivalence class of graphons. 
The same holds when $\widetilde{\Gamma}(w)$ is small,~i.e. there may be $u_1,u_2\in\blue{W_0}$ with $\delta_\Box(u_1,w)=\delta_\Box(u_2,w)=0$, $\widetilde{\Gamma}(w)=\Gamma(u_1)=\Gamma(u_2)$, but $\|u_1-u_2\|_\Box$ is large. This phenomenon is the root of complications in the proof of Theorem~\ref{thm-conj}, and the reason that a purely combinatorial proof could not be derived easily. 
\item[(ii)] A function $\Psi:W\rightarrow [0,1]$ is called a $\Gamma$-type function if it is $\|\cdot\|_\Box$-norm continuous, satisfies  the condition of 
Lemma~\ref{lem:Gamma-lower-bound}, and attains 0 when applied to Robinson graphons. It is very easy to verify that Theorem~\ref{thm:main} and Theorem~\ref{thm-conj} hold for any $\Gamma$-type function. 
In that sense, we think of $\Gamma$ as a prototype of $\Gamma$-type functions, whose definition is very natural. 
\end{itemize}
\end{remark}

\section{Rate of decay for $\widetilde {\Gamma}$ of samples from Robinson graphons}\label{sec:decay}
Let $w$ be a Robinson graphon, and consider a graph sequence $\{G_n\}_{n\in{\mathbb N}}$ with $G_n$ sampled from ${\mathbb G}(n,w)$. From Corollary 6.5 in  \cite{CGHJK}, we have that $\lim_{n\rightarrow \infty}\widetilde{\Gamma}(G_n)=0$ almost surely. In this section, we obtain upper and lower bounds for the speed of this convergence. 
In particular, we prove Theorem~\ref{thm:convergence}, which gives order bounds on $\widetilde{\Gamma}(G_n)$ for graph sequences sampled from  graphons with a flat region (rectangular region in $[0,1]^2$ on which $w$ is constant) and  steep graphons (graphons with partial derivatives bounded away from zero). Our results show that the decay is an order of magnitude faster for steep graphons than for graphons with a flat region. This confirms the intuitive notion that for graphons with stronger linearly embedded structures  $\widetilde {\Gamma}(G_n)$ converges to zero faster.

\subsection{Some definitions and notations for this section}
Due to the graph-theoretic nature of our results, we will replace the graphon-based definition $\widetilde {\Gamma} (G)$ by the graph-based parameter $\Gamma^*(G)$, originally defined in \cite{CGHJK} as follows:
Let $G$ be a graph with a linear order $\prec$ on its vertex set $V(G)$, and denote $|V(G)|=n$. For every $v\in V(G)$, the collection of all the neighbours of $v$ is denoted by $N(v)$. Also, the down-set  $D(v)$ and the up-set $U(v)$  of $v$ are defined as follows:
 $$D(v)=\{x\in V(G): x\prec v\}\ \mbox{ and }\ U(v)=\{x\in V(G):v\prec x\}.$$
\begin{definition}
Let $A\subseteq V(G)$, and $\prec$ be a linear order on the vertex set of $G$. We define,
\begin{eqnarray*}
\Gamma^* (G,\prec, A)&=&\frac{1}{n^3}\sum_{u\prec v}{\big[} 
|N(v)\cap A\cap D(u)|-|N(u)\cap A\cap D(u)|{\big]}_+\\
&+&\frac{1}{n^3}\sum_{u\prec v}{\big[}
|N(u)\cap A\cap U(v)|-|N(v)\cap A\cap U(v)|{\big ]}_+.
\end{eqnarray*}
We also define 
\[\Gamma^*(G,\prec)=\max_{A\subseteq V(G)}\Gamma^*(G,\prec,A)\ \mbox{ and } \ \Gamma^*(G)=\min_{\prec}\Gamma^*(G,\prec),\]
where the minimum is taken over all the linear orderings of $V(G)$.
\end{definition}
It follows from  \cite[Corollary 5.2]{CGHJK} that the parameters $\Gamma^*$ and $\widetilde{\Gamma}$ are asymptotically equal,~i.e.~\pink{for a sequence of graphs $\{ G_n\}$ of increasing order $n$},
\begin{equation}\label{Gamma*}
\Gamma^*(\pink{G_n})=\widetilde{\Gamma}(\pink{G_n})+O(\frac{1}{n}).
\end{equation}

We will use the term {\sl with exponential probability (w.e.p.)} to denote that an event holds asymptotically with probability at least $1-\exp{(-c\log^2 n)}$ for some positive constant $c$. Thus, a polynomial number of events that all hold {\sl w.e.p.\ }will hold simultaneously {\sl w.e.p.} as well. We will make extensive use of a well-known inequality, quoted below, which shows that the sum of bounded independent variables is close to its expected value.
\begin{theorem}[Hoeffding's inequality]\label{thm:Hoeffding-bounded}
Let $\{ X_i\}_{i=1}^N$ be a sequence of independent random variables bounded by the interval $[0,1]$, and let $X=\sum_{i=1}^N X_i$.  Then 
\begin{eqnarray}
\label{Hoeffding}
{\mathbb P} (|X-\sum_{i=1}^N \bbE(X_i)|\geq \alpha N)\leq 2 \exp{(-2\alpha^2 N)}.
\end{eqnarray}
\end{theorem}
Further concentration bounds used in this section can be found in \cite{Wainwright}. \blue{For the convenience of the reader, we also include these results in \ref{sec:appendix}.}

The proof of Theorem~\ref{thm:convergence} will follow from two propositions presented in the following two subsections, each of which addresses one case of the theorem.
\subsection{Graphons with a flat region}
In this section, we prove Proposition~\ref{prop:flat_graphon} from which case (i) of Theorem~\ref{thm:convergence} will follow. 
Namely, we give a lower bound for the  convergence rate of $\Gamma^*$ of graphs sampled from graphons with a constant rectangular region. 
In particular, this bound may be applied to the constant graphon, which is an extreme case of a graphon that is Robinson, but does not have a clear linear embedding.

Recall that $\Gamma^*$ is a graph parameter defined as $\min_{\prec}\Gamma^*(G,\prec)$,  where the minimum is taken over all the linear orderings of $V(G)$.
For any particular ordering $\prec$, one can use standard probabilistic techniques to show that \wep $\Gamma^*(G,\prec)\geq cn^{-1/2}$ for some constant $c$. 
This fact is due to the random fluctuation of the outcomes of the vertex pairs  whose labels fall in the constant region of $w$. 
In order to obtain a lower bound on $\Gamma^*(G)$, we have to establish that the same lower bound remains true for \emph{all} orderings. 
This cannot be achieved using a simple union bound argument, as the number of orderings of $V(G)$ is superexponential.
To address this issue, we partition the orderings into classes according to a very coarse partial order of the vertex set. Next, we  establish the desired lower bound on $\Gamma^*(G,\prec)$ for an entire class of orderings simultaneously, and then show that the number of classes is small enough so that a union bound argument guarantees that \wep the lower bound holds for all classes simultaneously.

We will need the following lemma.

\begin{lemma}
\label{lem:ordering}
Let $V_S,V_T$ be two subsets of a set $V$, and let $\prec $ be an ordering of $V$. Let $V_1\subseteq V_S$ and $V_2\subseteq V_T$ so that 
\begin{equation}\label{cond:initial segment}
 V_1\prec V_S\setminus V_1\text{ and } V_2\prec V_T\setminus V_2.
 \end{equation}
 Then
$V_1\prec V_T\setminus V_2$ or $V_2\prec V_S\setminus V_1$.
\end{lemma}
\begin{proof}
Let $V=\{ v_1,\dots ,v_n\}$ where the vertices are labeled according to $\prec$, i.e.~we have  $v_1\prec v_2\prec \dots \prec v_n$. 
Our assumption on $\prec$ implies that if $v_i\in V_1$ and $v_j\in V_S\setminus V_1$, then $i<j$; an equivalent statement holds for $V_T$ and $V_2$. 

For $i=1,2$, let $k_i$ be the largest index in $V_i$, that is, $k_i=\max\{ k\,:\,v_k\in V_i\}$.  
Suppose $k_1\leq k_2$. Then $V_1\subseteq \{ v_1,\dots ,v_{k_1}\}$, and by Condition (\ref{cond:initial segment}) we have $V_T\setminus V_2\subseteq \{ v_{k_2+1},\dots ,v_n\}$. Therefore $V_1\prec V_T\setminus V_2$. Similarly, if $k_2<k_1$ then $V_2\prec V_S\setminus V_1$.
\end{proof}

\begin{proposition}
\label{prop:flat_graphon}
Let $w$ be a graphon with a constant rectangular region, that is, there are measurable sets $S,T\subseteq [0,1]$ with positive measure and a constant $p\in (0,1)$ such that $w$ assumes $p$  a.e.~on $S\times T$. Then there exist constants $c,\alpha>0$ so that,
for a $w$-random outcome $G\sim \cG(n,w)$   we have
\[
\mathbb{P}(\Gamma^*(G)\geq c n^{-1/2})\geq 1- \exp (-\alpha n).
\]
Moreover, the constants $c,\alpha$ depend only on $p$ and $\min\{ |S|,|T|\}$.
\end{proposition}
\begin{proof}
Let $w$ be as stated, and denote $s:=\min\{ |S|,|T|\}$.  
Let $G\sim \cG(n,w)$ be an outcome with vertex set $V=\{ 1,2,\dots , n\}$ and labels  $x_1,\ldots, x_n$. Define 
\[
V_S=\{ i\in V\,:\, x_i\in S\}\text{ and }V_T=\{ i\in V\,:\, x_i\in T\}.
\]
Clearly, $w(x,y)=p$ for all $x\in V_S$, $y\in V_T$.
By Hoeffding's inequality, with probability at least $1-2\exp (-2(s/4)^2n)$, we have 
\begin{equation}\label{eq:v-size}
 |V_S|\geq \frac34 sn\  \mbox{ and  }\  |V_T|\geq \frac34 sn.
 \end{equation}
We assume this is the case.

Let $0< \epsilon < \frac{s}4 $. (An appropriate choice for $\epsilon$ will become apparent at the end of this proof.) 
Consider partitions of $V_S$ and $V_T$ into sets  $V_1^S,V_2^S,V_3^S$ and $V_1^T,V_2^T, V_3^T$ respectively, so that
\begin{equation}\label{cond:Vi}
|V_1^S|=|V_2^S|=|V_1^T|=|V_2^T|=\lfloor \epsilon n\rfloor\quad   \text{ and }\quad  (V_1^T\cup V_2^T)\cap V_3^S=\emptyset.
\end{equation}
 We will show that \wep for \emph{any} such partition, there exists a sufficiently large subset of $V_3^S$ containing vertices that have larger than expected number of neighbours in $V_1^T$ and smaller than expected number of neighbours in $V_2^T$, and vice versa. 

Let $m=\lfloor \epsilon n\rfloor$. For $i=1,2$, let  
\begin{eqnarray}\label{def:B_C}
B_i^S&=&\{ y\in V_{3}^S: |N(y)\cap V_{i}^T| \geq pm+\sqrt{m}\}\\
C_i^S&=&\{ y\in V_{3}^S: |N(y)\cap V_{i}^T| \leq pm-\sqrt{m}\}\notag 
\end{eqnarray}
We are interested in the sets $B_1^S\cap C_2^S$ and $B_2^S\cap C_1^S$. First, we establish a lower bound on their sizes.

 For each $y\in V_{3}^S$, the variable $|N(y)\cap V_{1}^T|$ is the result of $m $ Bernoulli trials, with probability of success equal to $p$, so ${\mathbb E} (|N(y)\cap V_{1}^T|)=p m $. Moreover, since $V_3^S\cap V_1^T=\emptyset$, these variables are independent. Thus by the central limit theorem, for  $y\in V_{3}^S$, the random variable $Z_{m}=\frac{|N(y)\cap V_1^T|- pm}{\sqrt{mp(1-p)}}$ converges in distribution to a standard normal random variable with {\it cdf} $\Phi$ as $m$, or equivalently $n$, converges to infinity. So \blue{if $m$ is large enough}, for every $y\in V_3^S$, 
 \begin{equation}\label{eq:prob-estim1}
 \mathbb{P}(y\in B_1^S)\geq \frac34 \Phi( \frac{1}{\sqrt{p(1-p)}})\ \mbox{ and }\ \mathbb{P}(y\in C_2^S)\geq \frac34 \Phi( \frac{1}{\sqrt{p(1-p)}}). 
 \end{equation}

Next, we show that \wep $|B_1^S\cap {C_2^S}|$ has  a suitable lower bound.
Let $a_p=\Phi( \frac{1}{\sqrt{p(1-p)}}).$
Since $V_1^{T}\cap V_2^{T}=\emptyset$, the events $\{ y\in B_1^S\}$ and $\{y\in C_2^S\}$ are independent, so  (\ref{eq:prob-estim1}) implies that $\mathbb{P}(y\in B_1^S\cap C_2^S)\geq (\frac34 a_p)^2$. 
Note that $|B_1^S\cap C_2^S|$ is a sum of Bernoulli random variables with probability of success at least $(\frac34 a_p)^2$, so $\mathbb{E}(|B_1^S\cap C_2^S|)\geq \frac{9}{16}a_p^2|V_3^S|$. 
From the condition $V_3^S\cap (V_1^T\cup V_2^T)=\emptyset$ in~\eqref{cond:Vi}, these Bernoulli random variables are independent as well, 
so by Hoeffding's inequality, we have that, with probability at least
 $1-2\exp(-2(a_p^2/16)^2|V_3^S|)$,  the size $|B_1^S\cap {C_2^S}|$ is bounded below by $(1/2)a_p^2|V_3^S|\geq (1/2)a_p^2(\frac34 s-2\epsilon)n$. 
Employing a similar argument for $|B_2^S\cap {C_1^S}|$ and using a union bound, we obtain that if  \eqref{eq:v-size} is satisfied, then with probability at least $1-4\exp (-2(a_p^2/16)^2 (\frac34 s-2\epsilon )n)$,  the set sizes  $|B_1^S\cap C_2^S|$ and   $|B_2^S\cap C_1^S|$ are bounded below by $(1/2)a_p^2(\frac34 s-2\epsilon)n$.

To show that these lower bounds hold \emph{w.e.p.}\ for \emph{every} choice of $V_i^S$, $V_i^T$, $i=1,2,3$, satisfying \eqref{cond:Vi}, we count the number of all such partitions. The number of ways to choose $V_1^S$ and $V_2^S$ of size $\lfloor\epsilon n\rfloor$ is ${\binom{|V_S|}{\lfloor \epsilon n\rfloor}} \binom{|V_S|-\lfloor\epsilon n\rfloor}{\lfloor \epsilon n\rfloor}$; recall that $V_3^S$ is determined by $V_1^S$ and $V_2^S$. 
It is known, as a consequence of Stirling's approximation of the factorial, that $\lim_{n\rightarrow \infty}\frac{1}{n}\log_2\binom{n}{m}=H(\epsilon)$, if $\lim_{n\rightarrow \infty} \frac{m}{n}=\epsilon$ (see \cite{spencer-book}). Here, $H$ is the binary entropy function defined as $H(\epsilon)=-\epsilon\log_2\epsilon-(1-\epsilon)\log_2(1-\epsilon)$.
Given that $H$ is increasing on $[0,1/2]$,
we obtain the following upper bound on the number of ways to choose $V_i^S,V_i^T$:
$$
{\binom{|V_S|}{\lfloor \epsilon n\rfloor}} \binom{|V_S|-\lfloor\epsilon n\rfloor}{\lfloor \epsilon n\rfloor}{\binom{|V_T|}{\lfloor \epsilon n\rfloor}} \binom{|V_T|-\lfloor\epsilon n\rfloor}{\lfloor \epsilon n\rfloor}
\leq {\binom{n}{\lfloor \epsilon n\rfloor}}^4\leq 2^{5H(\epsilon)n}
$$

Define
\[
f_p(x)=2 (a_p^{2}/16)^2 (\frac34s-2x )- 5\ln (2) {H(x )}.
\]
Note that $f_p$ is strictly decreasing on $[0,1/4)$, and has positive value at $x=0$.  
Let {$\epsilon^*$} be the value of $x$ in $[0,s/4)$ where $f$ attains zero, or $\epsilon^*=s/4$ if no such value exists. We have now established the following claim:
\begin{claim}\label{claim:intersection}
For each $0<\epsilon <\epsilon^*$, with probability at least $1-2\exp (-2(s/4)^2n)- 4\exp (-\blue{f_p}(\epsilon)n)$, 
$V_S$ and $V_T$ both have size at least $(3/4)sn$, and for \emph{every} choice of $V_1^S,V_2^S,V_1^T,V_2^T$ satisfying Condition \eqref{cond:Vi},  $|B_1^S\cap C_2^S|$ and $|B_2^S\cap C_1^S|$ are  bounded below by $(1/2)a_p^2(\frac34 s-2\epsilon)n$. 
\end{claim}

Now consider any ordering $\prec$ of $V$, and let $V_1^S$ be the first $m=\lfloor \epsilon n\rfloor$ elements of $V_S$ according to $\prec$, and $V_2^S$ the next $m$ elements of $V_S$. Similarly, let $V_1^T$ and $V_2^T$ be the first $m$ elements and the next $m$ elements of $V_T$.
By Lemma~\ref{lem:ordering}, it follows that $(V_{1}^S\cup V_2^S)\prec V_3^T$ or 
$(V_{1}^T\cup V_2^T)\prec V_3^S$. By switching $S$ and $T$ if necessary, we may assume \emph{wlog} that $(V_{1}^T\cup V_2^T)\prec V_3^S$. In particular, this implies that $(V_1^T\cup V_2^T)\cap V_3^S=\emptyset$, so Condition \eqref{cond:Vi} is satisfied.

We now have that, for $i=1,2$, 
\begin{eqnarray}
n^3\Gamma(G,\prec,V_i^T)&\geq& \sum_{x,y\in V_3^S,\,  x\prec y} \big[|N(y)\cap V_i^T|-|N(x)\cap V_i^T|\big]_+ \nonumber\\
&
\geq 
&
\sum_{x\in C_i^S,\,y\in B_i^S, \, x\prec y } \big[|N(y)\cap V_i^T|-|N(x)\cap V_i^T|\big]_+\nonumber\\
&\geq& 2\sqrt{m}\, \left|\left\{(x,y)\in C_i^S\times B_i^S:\ x\prec y\right\}\right|.\label{eq:newset}
\end{eqnarray}
Note that any pair $(x,y)$ with $x\in  B_2^S\cap C_1^S$ and $y\in B_1^S\cap C_2^S$ belongs either to the set described in \eqref{eq:newset} with $i=1$  (if $x\prec y$), or to the same set with $i=2$ (if $y\prec x$).
So 
$$n^3\Gamma(G,\prec,V_1^T)+n^3\Gamma(G,\prec,V_2^T)\geq  2\sqrt{m}\, |B_2^S\cap C_1^S|\,|B_1^S\cap C_2^S|,$$
which implies that $2\Gamma(G,\prec)\geq  \Gamma(G,\prec,V_1^T)+\Gamma(G,\prec,V_2^T)\geq \frac{2\sqrt{m}}{n^3}\, |B_2^S\cap C_1^S|\,|B_1^S\cap C_2^S|$.

Taking $\epsilon =\epsilon ^*/2$, $\alpha =(\ln 6)\left(\min\{\frac{s^2}{8},f(\epsilon^*/2)\}\right)$ and  $c=(1/5)a_p^4(\frac34 s-2\epsilon)^2 \sqrt{\epsilon}$, and using Claim~\ref{claim:intersection}, we get that 
with probability at least $1-\exp (-\alpha n)$, we have $\Gamma^*(G,\prec)\geq cn^{-1/2}$ for all orderings $\prec$.
\end{proof}

\subsection{Steep graphons}
The next proposition gives an upper bound on $\Gamma^*$ for graphs sampled from a steep graphon, and thus covers case (ii) of Theorem~\ref{thm:convergence}. Steep graphons are graphons for which the partial derivatives exist and are bounded away from zero.

Since we are establishing an upper bound on $\Gamma^*$, we only need to show the upper bound for $\Gamma^*(G,\prec)$ for some ordering $\prec$. Guided by  the Robinson property, we use the natural ordering  of vertices induced from the labels assigned by the random process $\cG(n,w)$. We formalize this concept as follows: For $G\sim \cG(n,w)$, we denote the vertex set of $G$ by $V=\{ 1,\dots, n\}$, and identify each vertex $i$ with the \emph{label} $x_i\in [0,1]$, which is its sampled value. We assume that the vertices are ordered so that $x_1\leq x_2\leq \dots \leq x_n$. 
The graph $G$, labeled in this manner, is called a \emph{labeled outcome} of $\cG(n,w)$. Since the labels assigned to different vertices are almost surely distinct, we can assume that the vertex ordering mentioned above is unique.

The rather complicated nature of our proof is due to the fact that $\Gamma^*(G,\prec)$ is defined as the maximum of $\Gamma^*(G,\prec , A)$ over all possible choices of $A$. Thus,  we must prove that the proposed upper bound dominates $\Gamma^*(G,\prec , A)$  for {\sl every}  subset $A$ of $V(G)$. To do so, we  proceed by partitioning the interval $[0,1]$ into a large number of small intervals. We then show that the following two facts hold with exponential probability: on the one hand, pairs of vertices chosen from \blue{intervals that are far apart} do not have a positive contribution to $\Gamma^*$, as $w$ is a steep Robinson graphon; on the other hand, the contribution to $\Gamma^*$ coming from pairs of vertices chosen from \blue{intervals that are close together} can be bounded from above efficiently.

\begin{proposition}
\label{prop:steep_graphon}
Let $w:[0,1]^2\rightarrow (0,1)$ be a Robinson graphon, whose partial derivatives on $[0,1]^2$ exist and are bounded away from 0. That is, 
there exists a positive constant $c$, such that  
\begin{equation}\label{lowerbound-c}
\Big|\frac{\partial w}{\partial x}\Big|, \Big|\frac{\partial w}{\partial y}\Big|\geq c.
\end{equation}
Then, for an outcome $G\sim \cG(n,w)$,
\[
\mathbb{P}\left(\Gamma^* (G)\leq \left(\frac{\pink{26}}{c}\right) n^{-2/3}\right)\geq 1-\exp (-\Omega (\blue{\log^2n})).
\]
\end{proposition}

\begin{proof}
Let  $m=\phi(n)$, where $\phi(n)$ is an integer-valued function so that $\lim_{n\rightarrow \infty}{\frac{\phi(n)}{\log^3(n)}>\blue{1}} $ and $\lim_{n\rightarrow \infty}\frac{\phi(n)\log n}{\sqrt{n}}=0$.
An appropriate choice for a function $\phi$, satisfying both of these properties, becomes apparent at the end of the proof. 
Assume $n$ is large enough so that 
\begin{equation}\label{cond:n}
\log^3(n)< \phi (n) < 0.1\sqrt{n}\log^{-1} n.
\end{equation}
Next, divide $(0,1]$ into $m$ equal-sized intervals; namely $I_j=(\frac{j}{m},\frac{j+1}{m}]$, for $0\leq j<m$.

Let $G\sim \cG(n,w)$ be a labeled outcome, i.e. $V(G)$ is labeled so that $0<x_1< x_2< \ldots<x_n<1$. 
Throughout the proof, we identify a vertex $i$ in $V$ with its label $x_i$ in $(0,1)$, which allows us to think of $V$ \blue{both as a set of vertices} and as a subset of $(0,1)$. 

The partition  $\{I_j\}_{0\leq j<m}$ of $(0,1]$ results in an analogous partition of the vertex set $V$, which we  denote by $\{V_j\}_{0\leq j<m}$. Namely, for every $0\leq j< m$, $V_j=\{i\in V: x_i\in I_j\}$. 
We will see that \emph{w.e.p.}, these sets will all be close to their expected size.
Precisely, fix $0\leq j<m$, and note that $|V_j|$ is the sum of $n$ independent Bernoulli variables with success probability $|I_j|$. 
Applying Hoeffding's inequality with parameter $\alpha=\frac{\log n}{\sqrt{n}}$, we conclude that, with probability at least $1-2\exp(-2\log^2n)$, we have 
\begin{equation}\label{eq:vertex-density}
\big| |V_j| - n|I_j|\big|\leq \sqrt{n}\log n.
\end{equation}
Since the number of intervals $I_j$ is sub-linear in $n$, we can assume that \emph{w.e.p.},\   (\ref{eq:vertex-density}) holds for every $0\leq j< m$. Because of the fact that $\sqrt{n}\log n\leq 0.1\frac{n}{\phi(n)}$ (from (\ref{cond:n})), together with $|I_j|=\frac{1}{\phi(n)}$, the inequality (\ref{eq:vertex-density}) implies that
$(0.9)\frac{n}{\phi(n)}\leq |V_j|\leq (1.1)\frac{n}{\phi(n)}$.
Thus, the event 
\begin{equation}\label{eq:vertex-density-Vi}
\mathcal{B}_1=\bigcap_{0\leq j<m} \Big\{ (0.9)\frac{n}{\phi(n)}\leq |V_j|\leq (1.1)\frac{n}{\phi(n)} \Big\}
\end{equation}
holds \emph{w.e.p.}

Recall that the definition of $\Gamma^*$ consists of two parts, which we call $\Gamma^*_\ell$ and $\Gamma^*_u$. We will only bound $\Gamma^*_\ell$; an identical bound applies to $\Gamma^*_u$. Given a set $A\subseteq V$, let $A_k= A\cap V_k$, for $k=0,\dots , m-1$, and view each $A_k$ as a subset of $[0,1]$. Note that $[0,x)\cap A_k=A_k$ when $x\in V_i$ and $k<i$. So we get, 
\begin{eqnarray}
n^3\Gamma^*_\ell (G,A,\prec) 
&=& \sum_{x,y\in V,\  x<y}  \left[ \sum_{k=0}^{m-1} |N(y)\cap [0,x)\cap A_k|-|N(x)\cap [0,x)\cap A_k|\right]_+ \nonumber\\
 &\leq & \sum_{x,y\in V,\  x<y}\,   \sum_{k=0}^{m-1}
 \, \Big[ |N(y)\cap [0,x)\cap A_k|-|N(x)\cap [0,x)\cap A_k|\Big]_+ \nonumber\\
&=& \sum_{0\leq k< i\leq j< m}\sum_{{\scriptsize{\begin{array}{c}
x\in V_i \\
y\in V_j\\
x<y\\
\end{array}}}} 
\Big[ |N(y)\cap A_k|-|N(x)\cap A_k| \Big]_+\label{eq:part1}\\
&+& \sum_{0\leq k= i\leq j< m}\sum_{{\scriptsize{\begin{array}{c}
x\in V_i \\
y\in V_j\\
x<y\\
\end{array}}}} \Big[ |N(y)\cap [0,x)\cap A_k|-|N(x)\cap  [0,x)\cap A_k| \Big]_+.\label{eq:part2}
\end{eqnarray}
To obtain an upper bound for the above sum, we need to prove a few concentration results on random variables of the form $|N(x)\cap S|$. The graph $G$ is determined by two sets of random variables: firstly the random variables $x_i$ with values in $[0,1]$ which are assigned as labels to the vertices of $G$, and secondly the binary variables $Y_{x,z}$, where $Y_{x,z}=1$ precisely when the pair of vertices labeled as $x$ and $z$ form an edge in $G$. According to the definition of $\cG(n,w)$,  the random variables $Y_{x,z}$ are independent Bernoulli variables with probability of success $w(x,z)$. 
Event $\mathcal{B}_1$ defined in (\ref{eq:vertex-density-Vi}) is only a function of the labels. From now on, we assume that this event occurs.

Recall that $V$ is ordered so that the labeling $x_1,\ldots,x_n$ is increasing.  For a vertex $x\in V$ and a subset $S\subseteq V\setminus\{x\}$, define the following random variable.
\begin{equation*}\label{deltaHoeffding}
\delta_{x_1,\ldots,x_n} (x, S)= \sum_{s\in S} Y_{s,x}-\sum_{s\in S} \bbE\, Y_{s,x}=|N(x)\cap S|-\bbE(|N(x)\cap S|).
\end{equation*}
%
 Note that $\bbE(|N(x)\cap S|)=\sum_{z\in S}w(x,z)$. We will now show how the definition of $\delta_{x_1,\ldots,x_n}$ allows us to bound $n^3\Gamma^*_\ell(G, A,\prec)$ further. To simplify notation, we denote $\delta_{x_1,\ldots,x_n} (x, S)$ by $\delta(x, S)$, when the assignment $\{x_i\}_{i=1}^n$ is understood.

Suppose $x\in V_i$, $y\in V_j$, $x<y$, and $k<i \leq  j$. Since $\bbE|N(x)\cap A_k|-\bbE|N(y)\cap A_k|=\sum_{z\in A_k}w(x,z)-w(y,z)$, we have
\begin{eqnarray*}
 \Big[ |N(y)\cap A_k|-|N(x)\cap A_k| \Big]_+ = \Big[ \delta(y,A_k)-\delta(x,A_k)-\sum_{z\in A_k}(w(x,z)-w(y,z))\Big]_+.
\end{eqnarray*}
Given that  $z<x<y$ and $w$ is Robinson, every summand $w(x,z)-w(y,z)$ in the above sum is non-negative.
Moreover, by the Mean Value Theorem and condition (\ref{lowerbound-c}), we have $w(x,z)-w(y,z)\geq c(y-x)$ for every $z\in A_k$. Thus, 
\begin{eqnarray}
 \Big[ |N(y)\cap A_k|-|N(x)\cap A_k| \Big]_+&=& \Big[ \delta(y,A_k)-\delta(x,A_k)-\sum_{z\in A_k}(w(x,z)-w(y,z))\Big]_+\nonumber\\
&\leq&\Big[ \delta(y,A_k)-\delta(x,A_k)-c|A_k|(y-x)\Big]_+\nonumber\\
&\leq& \Big[|\delta(y,A_k)|-\frac{c|A_k|(y-x)}{2}\Big]_+\nonumber \\
&&\hspace*{2cm} +\Big[|\delta(x,A_k)|-\frac{c|A_k|(y-x)}{2}\Big]_+.\label{eq-deltaN}
\end{eqnarray}

Since $\bbE (\delta (x,S))=0$ for each set $S$, an upper bound will depend on concentration results for sums of random variables of type $\delta(x,S)$ or $\delta^2 (x,S)$. The following three claims will establish the necessary bounds.
%
\begin{claim}\label{claim-w.e.p}
For a labeled outcome $G\sim \mathbb{G}(n,w)$ with labeling  $x_1,\ldots,x_n$ assigned to $V(G)$, 
define the  event
\begin{equation}\label{def:deltaSumGood}
\mathcal{B}_2=\bigcap_{0\leq k< j <m}\bigcap_{S\subseteq V_k}\Big\{\sum_{x\in V_j} \delta_{x_1,\ldots,x_n}^2(x,S)<  \blue{7} \frac{n}{\phi(n)}|S|\Big\}.
\end{equation}
For $n$ large enough so that (\ref{cond:n}) holds, 
we have 
\[
{\mathbb P}({\mathcal B}_2\,|\,\mathcal{B}_1)\geq 1- \exp\left(-\frac{0.0025\, n}{\phi(n)}\right).
\]
That is, assuming ${\mathcal B}_1$, 
\wep for every  $0\leq k<j <m$ and subset $S\subseteq V_k$, $\sum_{x\in V_j} \delta^2(x,S)<  \blue{7} \frac{n}{\phi(n)}|S|$.  
\end{claim}
\begin{proof}[Proof of claim]
Fix an  outcome $G\sim \mathbb{G}(n,w)$ with labeling $x_1,\ldots,x_n$ for $V$. 
Fix $k<j$, a subset  $S\subset V_k$ and a vertex $x\in V_j$.
Since the random variables $Y_{s,x}$, $s\in S$, are independent, by Hoeffding's inequality (\ref{Hoeffding}), we have for all $t>0$,
$${\mathbb P}(\delta(x,S)\geq t)\leq \exp({-2t^2}/{|S|}) \ \mbox{ and } \ {\mathbb P}(\delta(x,S)\leq -t)\leq \exp({-2t^2}/{|S|}).$$
This means that $\delta(x,S)$ \blue{satisfies the conditions of the proposition in~\ref{prop:app-sub-G} with $\sigma^2=\frac{{|S|}}{4}$, and $\delta^2(x,S)-\bbE(\delta^2(x,S)\ch)$ is a sub-exponential random variable with parameters $\nu=2\sqrt{2}|S|, \alpha=2|S|$.} Moreover, 
 $\{\delta^2(x,S)-\bbE(\delta^2(x,S)\ch)\}_{x\in V_j}$ are independent random variables, as $V_j\cap S=\emptyset$. Thus, by 
 \blue{the tail bound on sums of sub-exponential variables as in Theorem~\ref{thm:sub-exp-tail-bounds}, part (ii)}, we have
\begin{equation*}
{\mathbb P}\Big(\sum_{x\in V_j} \delta^2(x,S)-\sum_{x\in V_j} \bbE(\delta^2(x,S)\ch)<t |V_j|\ch\Big)\geq 1-\exp\big(-\frac{|V_j|}{2}\min\big\{\frac{t^2}{\blue{8}|S|^2}, \frac{t}{\blue{2}|S|}\big\}\big).
\end{equation*}
Taking $t=\beta|S|$ with $\beta\geq 4$, the above inequality implies that, with probability at least $1-\exp(-\frac{\beta}{\blue{4}}|V_j|)$, 
\begin{equation}\label{bound1}
\sum_{x\in V_j} \delta^2(x,S)-\sum_{x\in V_j} \bbE(\delta^2(x,S)\ch)<\beta |V_j|\,|S|.
\end{equation}

There are at most $m^2 2^{|V_k|}$ choices for $j,k$, and $S$.
Assuming $\mathcal{B}_1$ holds,  and using conditions \eqref{cond:n}  and \eqref{eq:vertex-density-Vi}, we have that 
$m^2\leq n \leq \exp(0.1\frac{\sqrt{n}}{\phi(n)})\leq\exp (\frac{1}{9}|V_k|) \leq \exp ((1-\ln 2)|V_k|)$. 
Using (\ref{eq:vertex-density-Vi}) again, we get:
$$
m^2 2^{|V_k|}\exp(-\blue{\frac{\beta}{4}}|V_j|)\leq  \exp(|V_k|)\exp(-\blue{\frac{\beta}{4}}|V_j|)\leq  \exp\left(-(
\blue{\frac{\beta}{4}}(0.9)-1.1)(\frac{n}{\phi(n)})\right)$$
is exponentially small if $\blue{\beta>\frac{4(1.1)}{0.9}}$. Taking \blue{$\beta=4.9$} satisfies this condition.

Finally, note that $\bbE(\delta^2(x,S)\ch)=\sum_{y\in S} {\rm Var}(Y_{x,y}\ch)=\sum_{y\in S}w(x,y)-w(x,y)^2\leq |S|$. Therefore, $\sum_{x\in V_j} \bbE(\delta^2(x,S)\ch)\leq  \frac{1.1n|S|}{\phi(n)}$, and this finishes the proof. 
\end{proof}

To bound summand (\ref{eq:part2}), we need concentration results for a new type of sets defined as follows. Given a set $A\subseteq V$, an index $k$ and a vertex $x\in V_k$,  let $A_k=A\cap V_k$ and $\blue{A_{k,x}}=A\cap V_k\cap [0,x)$.
\begin{claim}\label{claim-w.e.p2}
For a labeled outcome $G\sim \mathbb{G}(n,w)$ with labeling  $x_1,\ldots,x_n$ assigned to the vertex set $V$,
define the  event
\[
{\cal B}_3=\bigcap_{0\leq i <m}\ \bigcap_{T, A_i\subseteq V_i}\Big\{\blue{\big|}\sum_{x\in T}\delta_{x_1,\ldots,x_n}(x,\blue{A_{i,x}})\blue{\big|}< \left( \frac{n}{\phi(n)}\right)^{\frac{3}{2}}\Big\}.
\]
 For $n$ large enough so that (\ref{cond:n}) holds,  we have 
\[
{\mathbb P}({\cal B}_3\,|\,\mathcal{B}_1)\geq 1-\blue{2}\exp \left(-0.2\frac{n}{\phi (n)}\right).
\]
That is, assuming ${\cal B}_1$, \wep $\mathcal{B}_3$ holds. 
\end{claim}
\begin{proof}[Proof of claim]
Fix an index $0\leq i< m$, two sets $A_i,T\subseteq V_i$, and a labeling $x_1,\ldots,x_n$ for $V$. 
As before, we denote $\delta_{x_1,\ldots,x_n} (x, S)$ by $\delta(x, S)$.
Recall that 
\[
\sum_{x\in T}\delta(x,\blue{A_{i,x}})=\sum_{x\in T,z\in \blue{A_{i,x}}}(Y_{x,z}-w(x,z)),
\] 
where $Y_{x,z}$ denotes the indicator variable of $z\in N(x)$ with probability of success $w(x,z)$. 
A random variable $Y_{x,z}$ contributes to the above sum, only if $x\in T\subseteq V_i$, $z\in \blue{A_{i,x}}\subseteq V_i$ and $z<x$. Therefore, $\sum_{x\in T}\delta(x,\blue{A_{i,x}})$ is the sum of at most $|V_i|^2/2$ independent variables, each with expected value zero.
 Applying Hoeffding's inequality, we get 
 \begin{equation}\label{eq:new1}
{\mathbb P}\left(\blue{\big|}\sum_{x\in T}\delta(x,\blue{A_{i,x}})\blue{\big|}\geq t\ch\right)\leq \blue{2}\exp\left(\frac{-2t^2}{|V_i|^2/2}\right).
\end{equation}
Let $t=\frac34|V_i|^{3/2}$, so the above bound on the probability is exponentially small. Assume that $\mathcal{B}_1$ holds, and $n$ satisfies (\ref{cond:n}). Note that there are at most  $m 2^{2|V_i|}\leq \exp (2|V_i|)$ choices for $i$, $A_i$ and $T$. \blue{So, applying a union bound to \eqref{eq:new1}, we get} 
\[
{\mathbb P}\Big(\blue{\bigcup_{0\leq i <m}\ \bigcup_{T, A_i\subseteq V_i}}\Big\{\blue{\big|}\sum_{x\in T}\delta(x,\blue{A_{i,x}})\blue{\big|}\geq \frac34|V_i|^{3/2}\Big\} \ \Big|\,\blue{\mathcal{B}_1}\Big)\leq \blue{2} \exp (-\frac{9}{4}|V_i|+ 2|V_i|) \leq \blue{2} \exp (-0.2\frac{n}{\phi (n)}).
\]
Since $\frac34|V_i|^{3/2}\leq \frac34(1.1)^{3/2}\left( \frac{n}{\phi(n)}\right)^{3/2}<\left( \frac{n}{\phi (n)}\right)^{3/2}$, the result follows.
\end{proof}

Finally, we show that with exponential probability, certain sums of variables $\delta(y,S)$ can be bounded from below. 
\begin{claim}\label{claim-w.e.p3}
For a labeled outcome $G\sim \mathbb{G}(n,w)$ with labeling  $x_1,\ldots,x_n$ assigned to the vertex set $V$,
define the  event
\[
{\cal B}_4=\bigcap_{0\leq i <m}\bigcap_{S\subseteq V_i}\bigcap_{T\subseteq V\setminus S}\Big\{\sum_{y\in T}\delta_{x_1,\ldots,x_n}(y,S)\blue{< \frac{1.1n^{3/2}}{\sqrt{\phi(n)}}}\Big\}.
\]
For $n$ large enough so that (\ref{cond:n}) holds, we have
\[
{\mathbb P}({\cal B}_4\,|\,\mathcal{B}_1)\geq 1- \exp (-0.1 n).
\]
That is, assuming ${\cal B}_1$, \wep $\mathcal{B}_4$ holds.
\end{claim}
\begin{proof}[Proof of claim]
Fix $0\leq i <m$, $S\subseteq V_i$, $T\subseteq V\setminus S$, and  a labeling $x_1,\ldots,x_n$ for $V$. Let $\delta(y,S)$ denote $\delta_{x_1,\ldots,x_n}(y,S)$ for short. Similar to the proof of Claim~\ref{claim-w.e.p2}, we can write $\sum_{y\in T}\delta(y,S)$ as  a sum of $|T||S|$ independent random variables. 
Suppose $x_1,\ldots,x_n$ are so that $\mathcal{B}_1$ holds. Then, $|T||S|\leq \frac{1.1n^2}{\phi(n)}$, and  applying Hoeffding's inequality we get 
$${\mathbb P}\left(\sum_{y\in T}\delta(y,S)\blue{< \frac{1.1 n^{3/2}}{\sqrt{\phi(n)}}}\ch\right)\geq 1-\exp\left(\frac{-2.42\, n^3}{\phi(n)|T||S|}\right)\geq 1-\exp(-2.2n).$$ 
The result will then follow by applying a union bound, 
as there are at most $\phi(n)2^{\frac{1.1n}{\phi(n)}+n}$ choices for $i$, $S$ and $T$, and $\phi(n)2^{\frac{1.1n}{\phi(n)}+n}\exp(-2.2 n)\leq \exp(-0.1 n)$.
\end{proof}

We have now shown that the probability that $\mathcal{B}_1$ does not hold is exponentially small, and, if $\mathcal{B}_1$ holds, then the probability that any of $\mathcal{B}_2$, $\mathcal{B}_3$ or $\mathcal{B}_4$ does not hold is also exponentially small. Using a union bound, we can then conclude that \wep $\mathcal{B}_1$, $\mathcal{B}_2$, $\mathcal{B}_3$ and $\mathcal{B}_4$ all hold. We will assume in the rest of the proof that this is the case.

{\bf Bounding summand (\ref{eq:part1}):}
The assumption that $\mathcal{B}_1$ $\mathcal{B}_2$, $\mathcal{B}_3$ and $\mathcal{B}_4$ all hold allows us to bound 
$\Gamma^*_\ell(G,A,\prec)$ as required. First, we bound the summand (\ref{eq:part1}).
%
%
Suppose $x\in V_i$, $y\in V_j$, $x<y$, and $k<i \leq  j$. As we saw earlier in (\ref{eq-deltaN}),
\begin{eqnarray}
 \Big[ |N(y)\cap A_k|-|N(x)\cap A_k| \Big]_+
&\leq& \Big[|\delta(y,A_k)|-\frac{c|A_k|(y-x)}{2}\Big]_+\label{eq-contr1}\\
&+&\Big[|\delta(x,A_k)|-\frac{c|A_k|(y-x)}{2}\Big]_+\label{eq-contr2}.
\end{eqnarray}
Since $x\in V_i$ and $y\in V_j$, and $i\leq j$, we have 
\begin{equation}\label{eq:boundy-x}
\frac{j-i-1}{\phi(n)}\leq y-x\leq \frac{j+1-i}{\phi(n)}.
\end{equation}
Note that if $y-x\geq \frac{2|\delta(y,A_k)|}{c|A_k|}$ then summand (\ref{eq-contr1}) will attain zero. 
\pink{
Using the lower bound on $y-x$}, we can see that   summand (\ref{eq-contr1}) is nonzero only if $j-1-\frac{2\phi(n)|\delta(y,A_k)|}{c|A_k|}< i$. 
%
%
\pink{Fix $y\in V_j$. If 
$i\in \left(j-1-\frac{2\phi(n)|\delta(y,A_k)|}{c|A_k|},j-\frac{2\phi(n)|\delta(y,A_k)|}{c|A_k|}\right]$, 
the term $\Big[ |\delta(y,A_k)|-\frac{c|A_k|(y-x)}{2}\Big]_+$ in summand (\ref{eq-contr1}) is at most $\frac{c|A_k|}{2\phi(n)}$ for each $x\in V_i$. 
There is exactly one such value of $i$ and this value gives 
a total contribution of at most $\frac{c|V_i|\,|A_k|}{2\phi(n)}\leq \frac{c(1.1)^2n^2}{2\phi(n)^3}$ to \eqref{eq:part0-1}.  
}
\pink{On the other hand,} there are at most $\frac{2\phi(n)|\delta(y,A_k)|}{c|A_k|}$ values of $i$ such that \pink{$j-\frac{2\phi(n)|\delta(y,A_k)|}{c|A_k|}< i\leq j$.}

Similarly, \pink{if $y-x\geq \frac{2|\delta(x,A_k)|}{c|A_k|}$ then summand (\ref{eq-contr2}) will attain zero. Using the upper bound from \eqref{eq:boundy-x} we see that,} given $x\in V_i$, there are at most  $\frac{2\phi(n)|\delta({x},A_k)|}{c|A_k|}+\pink{1}$ choices for $j$ so that summand (\ref{eq-contr2}) can possibly have nonzero values for $y\in V_j$, \pink{and one of those values contributes at most $\frac{c(1.1)^2n^2}{2\phi(n)^3}$ to \eqref{eq:part0-1}.}

Finally, since $\Big[|\delta(y,A_k)|-\frac{c|A_k|(y-x)}{2}\Big]_+\leq |\delta(y,A_k)|$ \blue{is valid for every $x\leq y$,}
and $|V_i|\leq \frac{1.1\, n}{\phi(n)}$ for every $i$, we have

\begin{eqnarray}
\mbox{Summand } (\ref{eq:part1})
&\leq&\sum_{\scriptsize{\begin{array}{c}
0\leq i\leq j<m \\
x\in V_i,y\in V_j\\
x<y\\
k< i
\end{array}}}
 \Big[ |\delta(y,A_k)|-\frac{c|A_k|(y-x)}{2}\Big]_+
+
\Big[ |\delta(x,A_k)|-\frac{c|A_k|(y-x)}{2}\Big]_+\nonumber\\
&\leq&
2\sum_{\scriptsize{\begin{array}{c}
0\leq j<m \\
y\in V_j\\
k< j
\end{array}}}
\left(
\left(\frac{2\phi(n)|\delta(y,A_k)|}{c|A_k|}\right)\left(\frac{1.1\, n}{\phi(n)}\right)
|\delta(y,A_k)|
%
\pink{+
\frac{c(1.1)^2}{2}\left(\frac{n^2}{\phi(n)^3}\right)}\right)\label{eq:part0-1}\\
 &\leq&
\left( \frac{4.4\, n}{c}\sum_{\scriptsize{\begin{array}{c}
0\leq j<m \\
y\in V_j\\
k< j
\end{array}}}
\frac{\delta(y,A_k)^2}{|A_k|}\right) \pink{+ c(1.4)\left(\frac{n^3}{\phi(n)^2}\right)}. \label{eq:part1-1}
\end{eqnarray}
Using our assumption that $\mathcal{B}_2$ holds, and the fact that there are only $m=\phi(n)$ many choices for $j$ and $k$, we get 
\begin{equation}\label{eq:bound1}
\mbox{Summand } (\ref{eq:part1-1})\leq\phi(n)^2
\left(\frac{4.4\, n}{c}\right)\left(\frac{{7}\, n}{\phi(n)}\right)
\pink{+\left(\frac{c(1.4)n^3}{\phi(n)^2}\right)}=\left(\frac{30.8}{c} \right)\phi(n) n^2\pink{+c(1.4)\left(\frac{n^3}{\phi(n)^2}\right)}.
\end{equation}
%
{\bf Bounding summand (\ref{eq:part2}):}
Using our assumption that $\mathcal{B}_1$, $\mathcal{B}_3$ and $\mathcal{B}_4$ hold, we now bound summand (\ref{eq:part2}).
For $x\in V_i$, recall the notation $\blue{A_{i,x}}=[0,x)\cap A_i$. 
\blue{Similar to the argument leading to \eqref{eq-deltaN}, we have
\begin{eqnarray*}
    \Big[ |N(y)\cap A_{i,x}|-|N(x)\cap A_{i,x}|\Big]_+ &\leq& \Big[\delta(y,A_{i,x})-\delta(x,A_{i,x})-c|A_{i,x}|(y-x)\Big]_+.\\
\end{eqnarray*}}
\blue{Since $[\cdot]_+$ is sub-additive, we get}
\begin{eqnarray}
\mbox{Summand } (\ref{eq:part2})
&\leq&
\sum_{\scriptsize{\begin{array}{c}
0\leq i<m\\
x\in V_i\\
\end{array}}}
\sum_{\scriptsize{\begin{array}{c}
i\leq j<m \\
y\in V_j\\
x<y\\
\end{array}}}\Big[ \blue{-\delta(x,A_{\blue{i,x}})}-\frac{c|A_{\blue{i,x}}|(y-x)}{2}\Big]_+\label{part2-1}\\
&+&
\sum_{\scriptsize{\begin{array}{c}
0\leq i<m\\
x\in V_i\\
\end{array}}}
\sum_{\scriptsize{\begin{array}{c}
y\in V\\
x<y\\
\end{array}}}
 \Big[ \blue{\delta(y,A_{i,x})}-\frac{c|A_{\blue{i,x}}|(y-x)}{2}\Big]_+.\label{part2-2}
\end{eqnarray}
For a fixed $x\in V_i$, a nonzero contribution of $y\in V_j$ in summand (\ref{part2-1}) is possible only if $\delta(x,A_{\blue{i,x}})\blue{<}0$ and 
$\frac{j-i-1}{\phi(n)}< \frac{\blue{-}2\delta(x,A_{\blue{i,x}})}{c|A_{\blue{i,x}}|}$, as $\frac{j-i-1}{\phi(n)}$ is a lower bound for $y-x$. Therefore, there are at most \blue{$\frac{2\phi(n)|\delta(x,A_{\blue{i,x}})|}{c|A_{\blue{i,x}}|}$} choices for $j\geq i$ with nonzero contribution in  (\ref{part2-1}). 
Let $V_i^{\blue{-}}$ denote the set $\{x\in V_i: \delta(x,A_{\blue{i,x}})\blue{<}0\}$. 
Since for every $j$, $|V_j|$ is bounded above by $\frac{1.1\, n}{\phi(n)}$,  we have
\begin{eqnarray}
\text{Summand } (\ref{part2-1})&=& \sum_{\scriptsize{\begin{array}{c}
0\leq i<m\\
x\in V_i\\
\end{array}}}
\sum_{\scriptsize{\begin{array}{c}
i\leq j<m \\
y\in V_j\\
x<y\\
\end{array}}}\Big[\blue{-} \delta(x,A_{\blue{i,x}})-\frac{c|A_{\blue{i,x}}|(y-x)}{2}\Big]_+\nonumber\\
&\leq& \sum_{\scriptsize{\begin{array}{c}
0\leq i<m\\
x\in V_i^{\blue{-}}\\
\end{array}}}  \big(\frac{2\phi(n) \blue{|}\delta(x,A_{\blue{i,x}})\blue{|}}{c|A_{\blue{i,x}}|}\big)\big(\frac{1.1\, n}{\phi(n)}\big)\blue{|}\delta(x,A_{\blue{i,x}})\blue{|}\nonumber\\\
&=&\frac{2.2\, n}{c} \sum_{\scriptsize{\begin{array}{c}
0\leq i<m\\
x\in V_i^{\blue{-}}\\
\end{array}}} \frac{\delta(x,A_{\blue{i,x}})^2}{|A_{\blue{i,x}}|}.\nonumber
\end{eqnarray}
By definition, it is clear that $|\delta(x,A_{\blue{i,x}})|\leq |A_{i,x}|$. So,  we have 
\begin{eqnarray}\label{eq:bound2}
\text{Summand } (\ref{part2-1})&\leq&\blue{\frac{2.2\, n}{c} \sum_{0\leq i< m}\sum_{x\in V_i^{\blue{-}}} -\delta(x,A_{i,x})=
\frac{2.2\, n}{c} \sum_{0\leq i< m}|\sum_{x\in V_i^{\blue{-}}} \delta(x,A_{i,x})|}\nonumber \\
&\leq&\frac{2.2\, n}{c}\phi(n)\big(\frac{n^{\frac{3}{2}}}{\phi(n)^{\frac{3}{2}}}\big)=\frac{2.2}{c}\big(\frac{n^{\frac{5}{2}}}{\sqrt{\phi(n)}}\big),
\end{eqnarray}
\noindent where the last inequality comes from our assumption that $\mathcal{B}_{3}$ holds. This bounds summand (\ref{part2-1}).
To bound summand (\ref{part2-2}), we need to take a different approach, as the distance between $x$ and $y$ in this case, cannot be bounded only in terms of \blue{$y$}. So, we proceed as follows. For each $x$, define $T_x=\{y\in V: y>x\  \mbox{ and }\ \delta(y,A_{\blue{i,x}})\blue{>}0\}.$ Note that
\begin{eqnarray}
\text{Summand } (\ref{part2-2})&=&\sum_{\scriptsize{\begin{array}{c}
0\leq i<m\\
x\in V_i\\
\end{array}}}
\sum_{\scriptsize{\begin{array}{c}
y\in V\\
x<y\\
\end{array}}}\Big[\blue{\delta(y,A_{i,x})}-\frac{c|A_{i,x}|(y-x)}{2}\Big]_+\nonumber\\
&\leq& \sum_{\scriptsize{\begin{array}{c}
0\leq i<m\\
x\in V_i\\
\end{array}}}
\sum_{y\in T_x}\blue{\delta(y,A_{i,x})}.\nonumber
\end{eqnarray}
Given the assumption that $\mathcal{B}_{4}$ from Claim~\ref{claim-w.e.p3} holds, we have
\begin{eqnarray}\label{eq:bound3}
\text{Summand } (\ref{part2-2}) \leq \sum_{\scriptsize{\begin{array}{c}
0\leq i\leq m-1\\
x\in V_i\\
\end{array}}}
\sum_{y\in T_x}\blue{\delta(y,A_{i,x})}
\leq n\left(\frac{{1.1}n^{\frac{3}{2}}}{\sqrt{\phi(n)}}\right)=\frac{{1.1}n^{\frac{5}{2}}}{\sqrt{\phi(n)}}.
\end{eqnarray}
Putting inequalities (\ref{eq:bound1}), (\ref{eq:bound2}) and (\ref{eq:bound3}) together, we conclude that 
$$
n^3\Gamma^*_\ell(G, A,\prec)\leq \left(\frac{30.8}{c} \right)\phi(n) n^2+\pink{c(1.4)\left(\frac{n^3}{\phi(n)^2}\right)}+\frac{{2.2}}{c}\left(\frac{n^{\frac{5}{2}}}{\sqrt{\phi(n)}}\right)+\frac{{1.1}\, n^{\frac{5}{2}}}{\sqrt{\phi(n)}}.
$$
Thus, the best upper bound, in order, will be obtained when $\phi(n)=\theta(n^{1/3})$, so that $\phi(n) n^2$, \pink{$n^3/\phi(n)^2$} and $\frac{n^{\frac{5}{2}}}{\sqrt{\phi(n)}}$ are of the same order.
From the conditions on $w$ and the fact that $w\geq 0$, we have that $c\leq 1$. Using this, and taking $\phi (n)=\pink{\frac{1}{2}}n^{\frac{1}{3}}$,
we have 
$$n^3\Gamma^*_\ell(G, A,\prec)\leq \left(\frac{30.8}{\pink{2}}\pink{+4(1.4)}+\pink{\sqrt{2}}(2.2+1.1)\right)\left(\frac 1{c}\right) n^{7/3}\leq \left(
\frac{\pink{26}}{c}\right)n^{7/3}.$$
\end{proof}

\pink{
\section{Conclusions and further work}

We have established that the graphon parameter $\Gamma$, first defined in \cite{CGHJK}, is indeed a suitable gauge for measuring the extent to which a graphon lacks the Robinson property.
It was already known that the function $\Gamma$ attains zero precisely when applied to a Robinson graphon. 
In this paper, we show that for any graphon $w$, there is a Robinson graphon $R_w$ so that $\|R_w-w\|_\Box\leq 14\Gamma(w)^{1/7}$. Thus, any graphon with small $\Gamma$-value is close, in cut norm, to a Robinson graphon. As a corollary of this approximation result, we prove that $\Gamma$ detects graph sequences sampled from Robinson graphons.  That is, for a convergent graph sequence $\{G_n\}_{n\in {\mathbb N}}$, we have $\Gamma(G_n)\to 0$ \emph{iff} the limit object has a Robinson representative. In this paper, we have focused on proving the existence of the Robinson approximation $R_w$ and the corresponding bound $\Gamma(w)^{1/7}$. We believe that a worthwhile direction of further study is to improve the exponent on the bound.

Computing, or even approximating, the Robinson approximation $R_w$ is another important topic of investigation, which is beyond the scope of this work. The definition of $R_w$ involves calculating the supremum of the average of $w$ over sets of certain shape and size. Even for graphons corresponding to matrices, it is not clear to us how $R_w$ can be efficiently computed. We believe there may be efficient algorithms to approximate $R_w$, and we find this topic worthy of further study.

Our results show that for a converging sequence of graphs, the equivalence class of the limiting graphon  contains a Robinson graphon \emph{iff} the sequence of $\Gamma$-values of the graphs converge to zero. However, to compute $\Gamma$ of a sampled graph, one would have to come up with a labeling of the graph that approximates the Robinson ordering of the graphon. This relates to the problem of {\sl graph seriation}, with important practical applications. We hope to use the theory developed in this paper as a stepping stone for making progress on the seriation problem for Robinson similarity data with errors. 

Finally, we show that the convergence rates of $\Gamma$ for sequences of random graphs sampled from a Robinson graphon depend on how strongly the graphon exhibits the Robinson property. The relevant result for steep graphons, only applies to graphons that are nowhere zero. Graphons that are zero outside a band around the main diagonal form a natural model for many applications. It would be worthwhile to extend this result to such graphons. 
}
\section*{Acknowledgements}
The first author acknowledges support from National Science Foundation Grant DMS-1902301 during the preparation of this article.
The second author was funded by an NSERC Discovery Grant and a sabbatical grant from Dalhousie University.
The authors completed this project while on a Research in Teams visit at the Banff International Research Center. We are grateful to BIRS for the  financial support and hospitality. 
\blue{Finally, we sincerely thank the anonymous reviewers for critically reading the manuscript and suggesting numerous improvements.}

\appendix
\blue{
\section{}\label{sec:appendix}
This appendix includes some of the definitions and concentration results that we need for Section~\ref{sec:decay}. A beautiful exposition of this material can be found in \cite[Section 2.1]{Wainwright}.

Let $X$ be a Gaussian random variable with mean $\mu$ and variance $\sigma^2$. It is easy to see that $X$ has the moment generating function ${\mathbb E}(e^{\lambda X})=e^{\mu\lambda+\frac{\sigma^2\lambda^2}{2}}\ \mbox{ valid for all } \lambda\in {\mathbb R}.$
The following upper deviation inequality for $X$ follows directly from Chernoff bound:
\begin{equation*}
{\mathbb P}(X\geq \mu+t)\leq e^{-\frac{t^2}{2\sigma^2}}, \mbox{ for all } t\geq 0.
\end{equation*}
This deviation/concentration result can be generalized to a large class of non-Gaussian random variables, which mimic the behaviour of Gaussian random variables to some extend; such random variables are described in the following definition.
\begin{definition}
A random variable $X$ with mean $\mu$ is said to be \emph{sub-Gaussian} if there exists $\sigma>0$ such that for all $\lambda\in {\mathbb R}$, we have
\begin{equation}
{\mathbb E}(e^{\lambda(X-\mu)})\leq e^{\frac{\sigma^2\lambda^2}{2}}.
\end{equation}
The constant $\sigma$ satisfying the above equation is called the \emph{sub-Gaussian parameter}, and $\sigma^2$ is called the \emph{variance proxy}. 
\end{definition}
In this article, we use equivalent characterizations of sub-Gaussian random variables as listed below.
\begin{theorem}[Equivalent characterizations of sub-Gaussian random variables](Cf.~\cite[Theorem 2.6]{Wainwright})
\label{thm:equivalence-sub-Gaussian}
Let $X$ be a random variable with zero mean. Then the following are equivalent.
\begin{itemize}
    \item[(i)] $X$ is a sub-Gaussian random variable with parameter $\sigma\geq 0$, that is,  for all $\lambda\in {\mathbb R}$ we have
    $${\mathbb E}(e^{\lambda X})\leq e^{\frac{\sigma^2\lambda^2}{2}}.$$
    \item[(ii)] There is a constant $c\geq 0$ and Gaussian random variable $Z$ with zero mean and variance $\tau^2$ such that 
    $${\mathbb P}(|X|\geq s)\leq c\, {\mathbb P}(|Z|\geq s), \mbox{ for all } s\geq 0.$$
    \item[(iii)] There is a constant $\theta\geq 0$ such that 
    $${\mathbb E}(X^{2k})\leq \frac{(2k)!}{2^k k!}\theta^{2k}, \mbox{ for all } k=1, 2,\ldots .$$
\end{itemize}
\end{theorem}
Clearly, any Gaussian random variable with variance $\sigma^2$ is sub-Gaussian with parameter $\sigma$. More generally, any bounded random variable, supported in some interval $[a,b]$, is sub-Gaussian with parameter at most $\frac{b-a}{2}$ (\cite[Exercise 2.4]{Wainwright}). An important concentration bound for sums of independent sub-Gaussian random variables is provided by Hoeffding's inequality. A special case of this inequality is stated in Theorem~\ref{thm:Hoeffding-bounded}, where Hoeffding's bound is applied to the case of bounded random variables supported in $[0,1]$.
\begin{theorem}[Hoeffding bound]\cite[Proposition 2.5]{Wainwright}
Let $\{X_i\}_{i=1}^N$ be a collection of independent sub-Gaussian random variables with sub-Gaussian parameters $\sigma_i$. Then for all $t\geq 0$, we have
\begin{equation*}
    {\mathbb P}\left(\sum_{i=1}^n X_i-\sum_{i=1}^n{\mathbb E}(X_i)\geq t\right)\leq \exp\big(-\frac{t^2}{2\sum_{i=1}^n\sigma_i^2}\big).
\end{equation*}
An identical upper bound holds for the left-hand-side event $\big\{\sum_{i=1}^n X_i-\sum_{i=1}^n{\mathbb E}(X_i)\leq t\big\}$.
\end{theorem}
A relaxation of the notion of sub-Gaussian random variables leads to the following definition.
\begin{definition}
A random variable $X$ with mean $\mu$ is said to be \emph{sub-exponential} if there exist non-negative parameters $(\nu,\alpha)$ such that for all $\lambda\in (-\frac{1}{\alpha},\frac{1}{\alpha})$, we have
\begin{equation}
{\mathbb E}(e^{\lambda(X-\mu)})\leq e^{\frac{\nu^2\lambda^2}{2}}.
\end{equation}
\end{definition}
Clearly, any sub-Gaussian random variable is sub-exponential as well. However, many sub-exponential random variables are not sub-Gaussian. 
Some equivalent characterizations of sub-exponential random variables are listed in the following theorem.
\begin{theorem}[Equivalent characterizations of sub-exponential random variables](Cf.~\cite[Theorem 2.13]{Wainwright})
\label{thm:equivalence-sub-exp}
Let $X$ be a random variable with zero mean. Then the following are equivalent.
\begin{itemize}
    \item[(i)] $X$ is a sub-exponential random variable with parameters $(\nu,\alpha)$, that is,  for all $\lambda\in (-\frac{1}{\alpha},\frac{1}{\alpha})$, we have
$${\mathbb E}(e^{\lambda X})\leq e^{\frac{\nu^2\lambda^2}{2}}.$$
    \item[(ii)] There is a constant $c_0>0$ such that ${\mathbb E}(e^{\lambda X})<\infty$ for all $|\lambda|\leq c_0$.
    \item[(iii)] There are constants $c_1,c_2> 0$ such that 
    $${\mathbb P}(|X|\geq s)\leq c_1e^{-c_2s}, \mbox{ for all } s>0.$$
    \item[(iv)] The value  $\gamma=\sup_{k\geq 2} \left(\frac{{\mathbb E}(X^{k})}{k!}\right)^{\frac{1}{k}}$ is finite. 
\end{itemize}
\end{theorem}
Similar to sub-Gaussian random variables, sub-exponential random variables also satisfy certain deviation/concentration inequalities, as stated below.
\begin{theorem}[Tail bound for (sums of independent) sub-exponential random variables]  (Cf. \cite[Section 2.1.3]{Wainwright})
\label{thm:sub-exp-tail-bounds}
\begin{itemize}
\item[(i)]  Let $X$ be a sub-exponential random variable with parameters $(\nu,\alpha)$. Then
$${\mathbb P}\big(X-{\mathbb E}(X)\geq t\big)\leq\left\{\begin{matrix}
    e^{-\frac{t^2}{2\nu^2}} &  0\leq t\leq \frac{\nu^2}{\alpha}\\
    e^{-\frac{t}{2\alpha}} & t>\frac{\nu^2}{\alpha}
\end{matrix}\right. $$
Consequently, for every $t>0$, we have ${\mathbb P}\big(X-{\mathbb E}(X)\geq t\big)\leq e^{-\frac{1}{2}\min(\frac{t}{\alpha},\frac{t^2}{\nu^2})}$.
An identical upper bound holds for the corresponding left-hand-side event.
\item[(ii)] Let $X_1,\ldots,X_n$ be independent random variables, such that each $X_i$ is sub-exponential with parameters $(\nu_i,\alpha_i)$. Then $\sum_{i=1}^n(X_i-{\mathbb E}(X_i))$ is sub-exponential with parameters
$\nu_*:=\sqrt{\sum_{i=1}^n\nu_i^2}$ and $\alpha_*:=\max_{1\leq i\leq n}\alpha_i$. Moreover, we have 
\begin{equation*}\label{eq:concentration-sum of-sub-exp}
{\mathbb P}\left(\frac{1}{n}\sum_{i=1}^n(X_i-{\mathbb E}(X_i))\geq t\right)
\leq\left\{\begin{matrix}
    e^{-\frac{nt^2}{2{\nu_*^2/n}}} &  0\leq t\leq \frac{\nu_*^2}{n\alpha_*}\\
    e^{-\frac{nt}{2\alpha_*}} & t>\frac{\nu_*^2}{n\alpha_*}
\end{matrix}\right. 
\end{equation*}
Consequently, for every $t>0$, we have 
${\mathbb P}\left(\frac{1}{n}\sum_{i=1}^n(X_i-{\mathbb E}(X_i))\geq t\right)\leq e^{-\frac{n}{2}\min(\frac{t}{\alpha_*},\frac{t^2}{{\nu_*^2/n}})}$.
An identical upper bound holds for the corresponding left-hand-side event.
\end{itemize}
\end{theorem}

In this article, the notions of sub-Gaussian and sub-exponential random variables and the relation between them are used, in particular, in the following format.
\begin{proposition}\label{prop:app-sub-G}
Let $X$ be a random variable with zero mean. Suppose there exists $\sigma>0$ such that for all $s>0$, we have 
$${\mathbb P}(X\geq s)\leq \exp\big(-s^2/2\sigma^2\big) \ \mbox{ and } \ {\mathbb P}(X\leq -s)\leq \exp\big(-s^2/2\sigma^2\big).$$
Then $X$ is a sub-Gaussian random variable with parameter $2\sqrt{2}\sigma$. Moreover, the random variable $X^2-{\mathbb E}(X^2)$ is a sub-exponential random variable with parameters $(8\sqrt{2}\sigma^2, 8\sigma^2)$.
\end{proposition}
\begin{proof}
We observe that ${\mathbb P}(|X|\geq s)\leq 2\, \exp\big(-s^2/2\sigma^2\big)$, thus, we have
\begin{eqnarray*}
{\mathbb E}(X^{2k})&=&\int_0^\infty {\mathbb P}(X^{2k}>s)\, ds =\int_0^\infty {\mathbb P}(|X|>s^{\frac{1}{2k}})\, ds\\
&\leq& 2\int_0^\infty \exp\big(-s^{\frac{1}{k}}/2\sigma^2\big)\, ds
= 2^{k+1}k\sigma^{2k}\int_0^\infty t^{k-1}\exp(-t)\, dt=2^{k+1}\sigma^{2k} k!,
\end{eqnarray*}
where we have used the change of variable $s^{\frac{1}{k}}/2\sigma^2=t$ in the penultimate equality. So, 
\begin{equation}\label{eq:bound-for-moments}
{\mathbb E}(X^{2k})\leq \frac{(2\sigma)^{2k} k!}{2^{k-1}},
\end{equation}
which implies that (iii) of Theorem~\ref{thm:equivalence-sub-Gaussian} holds, when $\theta=2\sigma$. Now from the proof of implication (iii) $\Rightarrow$ (i) of Theorem~\ref{thm:equivalence-sub-Gaussian} (see e.g.~Page 46 of \cite{Wainwright}), we conclude that 
$X$ is a sub-Gaussian random variable with parameter at most $2\sqrt{2}\sigma$.

To prove the second statement, let $Z=X^2-{\mathbb E}(X^2)$. Note that for every $k\in {\mathbb N}$, convexity of  $f(x)=x^k$, together with the triangle inequality, implies that $|X^2-{\mathbb E}(X^2)|^k\leq 2^{k-1}(X^{2k}+({\mathbb E}(X^2))^k)$. Moreover, Jensen's inequality guarantee that $({\mathbb E}(X^2))^k\leq {\mathbb E}(X^{2k})$. So, we have
\begin{eqnarray*}
{\mathbb E}(e^{\lambda Z})&\leq&{\mathbb E}(e^{|\lambda Z|})=1+\sum_{k=2}^\infty \frac{|\lambda|^k}{k!}{\mathbb E}\big|X^2-{\mathbb E}(X^2)\big|^k
\leq1+\sum_{k=2}^\infty \frac{|\lambda|^k 2^k}{k!}{\mathbb E}(X^{2k})\\
&\leq& 1+\sum_{k=2}^\infty \frac{|\lambda|^k 2^k}{k!}\left(2^{k+1}\sigma^{2k} k!\right)=
1+2(4|\lambda|\sigma^2)^2\sum_{k=0}^\infty (4|\lambda|\sigma^2)^k,\\
\end{eqnarray*}
where we have used the bound in \eqref{eq:bound-for-moments} for ${\mathbb E}(X^{2k})$. Since the sum $\sum_{k=0}^\infty (4|\lambda|\sigma^2)^k$ is bounded by 2, when $4|\lambda|\sigma^2\leq \frac{1}{2}$, we get 
$${\mathbb E}(e^{\lambda Z})\leq 1+4(4|\lambda|\sigma^2)^2\leq e^{64\sigma^4\lambda^2}, \ \mbox{ whenever } |\lambda|\leq \frac{1}{8\sigma^2},$$
that is, $Z$ is sub-exponential with parameters $(8\sqrt{2}\sigma^2, 8\sigma^2)$.
\end{proof}
}
%

\newcommand{\etalchar}[1]{$^{#1}$}

\end{document}